# A random observation-based management model of population dynamics and its ecological application

Hidekazu Yoshioka[1,2 *], Yuta Yaegashi[3], and Motoh Tsujimura[4]

[1] Graduate School of Natural Science and Technology, Shimane University, Nishikawatsu-cho 1060, Matsue, 690-8504, Japan

[2] Fisheries Ecosystem Project Center, Shimane University, Nishikawatsu-cho 1060, Matsue, 690-8504, Japan

[3] Graduate School of Agriculture, Kyoto University, Kitashirakawa-Oiwake-cho, Kyoto, 606-8502, Japan

[4] Graduate School of Commerce, Doshisha University, Karasuma-Higashi-iru, Imadegawa-dori, Kamigyo-ku, Kyoto, 602-8580, Japan

* Corresponding author: yoshih@life.shimane-u.ac.jp

**Abstract**

A new stochastic control problem of population dynamics under partial observation is formulated and analyzed both mathematically and numerically, with an emphasis on environmental and ecological problems. The decision-maker can only randomly and time-discretely observe and impulsively intervene the population dynamics governed by a regime-switching stochastic differential equation. The hybrid nature of the problem leads to an optimality equation containing an integro-differential equation and a static optimization problem. It is therefore different from the conventional Hamilton-Jacobi-Bellman equations. Existence and solvability issues of this optimality equation are analyzed in a viscosity sense. Its exact solution to a reduced but still nontrivial model is derived as well. The model is finally applied to a realistic environmental management problem in a river using a finite difference scheme.

## 1. Introduction

### 1.1 Problem background

Managing population dynamics of biotic communities is a central topic in various environmental and ecological studies. Biological resources management is one of such topics closely linked to human lives. It involves planning not only bio-economic harvesting policies of the resources but also cost-effective countermeasures against the predators (Ang et al., 2019; Blubaugh et al., 2016; Rindorf et al., 2016). Habitat environmental conditions of the biological resources are critical factors for successful bio-economic resource management as well. As an example of fisheries resources management, food availability (Huang et al., 2018; Steffensen and Mestl, 2016) and water quality (Asatryan et al., 2019; Di et al., 2019) should be monitored because of their direct influences on growth and survival of the resources.

The focus of this paper is a specific problem related to inland fisheries resources management. The problem is the nuisance algae bloom management in dam-downstream rivers: an urgent environmental problem in many rivers (Cullis et al, 2012; Yoshioka, 2019). Human interventions to river environment through dam operations significantly modify the river flow regime and often lead to a flow regime having smaller amount and magnitude than the natural ones. Modified flow regimes often provide habitable environment for filamentous green algae and invasive diatoms, such as *Cladophora glomerata* and *Didymosphenia Geminata* (Gladyshev and Gubelit, 2019; Gillis et al., 2018). Their bloom leads to ecological issues including habitat invasion against aquatic plants (Reaver et al., 2019), disturbance on biogenic compound and nutrient cycles in aquatic ecosystems (Picosz et al., 2017), and fisheries impacts with economic losses (Beville et al., 2012). River restoration including management of the algae bloom is therefore an important task for sustainable water environmental and fisheries management.

### 1.2 Mathematical background

Stochastic optimal control (Øksendal and Sulem, 2019) is a central mathematical framework for management of a broad class of system dynamics such as population dynamics. In this framework, the dynamics to be controlled are usually described with some stochastic differential equations (SDEs) (Yin and Zhu, 2009; Øksendal and Sulem, 2019).

In the standard cases, modeling and analysis of such control problems assume that the decision-maker, the controller of the system dynamics, can obtain a continuous information flow of the dynamics. Such examples include but are not limited to management of predatory and invasive species (Marten and Moore, 2011; Yaegashi et al., 2018), harvesting resources (Kvamsdal et al., 2016; Pinheiro et al., 2016), and environmental management (Bretschger and Vinogradova, 2019; Jiang et al., 2019). However, the assumption of continuously obtaining the information is not satisfied in many applications. This issue is critical especially if we consider population dynamics and/or environmental indices in natural environment because of technical difficulties to collect data and information processing and labor and instrumental costs as implied in field survey studies (Bolius et al., 2019; Cormier et al., 2018; Gable et al., 2017; Hopper et al., 2017; Ren et al., 2016). A more realistic model should assume that the system evolves continuously in time while the decision-maker can observe its dynamics only discretely in time (Yoshioka and Tsujimura, 2020). This assumption naturally leads to a control problem based on discrete observation as a kind of partial observations. The hybrid nature of the problem based on the continuous-time system dynamics and discrete-time observation leads to a more complicated model than the standard ones, but can be still handled with the dynamic programming principle.

Dyrssen and Ekström (2018) and Wang (2018) considered optimal stopping problems of diffusion-processes focusing on statistical inference and option pricing, in which a decision-maker decides the observation times based on the available information. Their framework has later been extended to fisheries (Yoshioka et al., 2019) and ecological management problems where a decision-maker can make interventions at the observation times (Yoshioka and Tsujimura, 2020). Another framework assumes exogenous random observation times. The randomness in the context of environmental management would come from difficulties to exactly specify the timing of interventions by busy environmental managers due to unavailability of workforce. Pham and Tankov (2008, 2009) carried out mathematical formulation and analysis, and computation of an optimal consumption problem based on the random observations. Their approach has later been extended to more complex problems in economics involving consumption and investment of multiple assets (Federico and Gassiat, 2014; Federico et al., 2017) and regime-switching markets (Gassiat et al., 2014). Both frameworks are more realistic than those with continuous observations. Focusing on epidemiological problems, Winkelmann et al. (2014) and Duwal et al. (2015) considered Markov decision processes under partial observations, and formulated analogous stochastic control problems. There exist related mathematical models of stochastic optimal control subject to Poisson intervention times (Wang, 2001; Castellano and Cerqueti, 2012; Lempa et al., 2014; Strini and Thonhauser, 2019). However, such approaches are still germinating in environmental and ecological study areas. Especially, problems with the random observations have been least explored to the best of the authorsôknowledge. This is our main motivation.

A common mathematical difficulty in analyzing the above-mentioned stochastic optimal control problems based on partial observation is that solving a control problem cannot be reduced to a classical Hamilton-Jacobi-Bellman: a degenerate elliptic or parabolic equation (Øksendal and Sulem, 2019). This is due to partial availability of information on the dynamics leading to a non-standard dynamic programming principle (Dyrssen and Ekström, 2018; Pham and Tankov, 2009). The optimality equations in these problems present mixed forms containing differential (or integro-differential) equations and static optimization problems having more complicated forms than the standard ones. Nevertheless, feasible numerical methods based on finite difference schemes and iterative methods are available (Yoshioka and Tsujimura, 2020; Pham and Tankov, 2008). Existence and solvability of the non-standard optimality equations can be addressed based on the concept of viscosity solutions (Crandall et al., 1992) with suitable modifications. These mathematical backgrounds motivate us to formulate and analyze ecological and environmental management problems under partial observations.

### 1.3 Contributions of this paper

The objectives of this paper are formulation, analysis, and application of a management problem of population dynamics under random observations. We focus on management of filamentous algae population dynamics in a dam-downstream river environment: a common environmental issues in many rivers worldwide (Cullis et al, 2012; Yoshioka, 2019). The system dynamics are described as an SDE driven by a continuous-time Markov chain representing a river flow process. The system is therefore piecewise deterministic (Yin and Zhu, 2009). The decision-maker, the manager of the river environment, can observe the population dynamics randomly following a controlled Poisson process, and can harvest (reduce) the population at each observation time. Population dynamics having regime-switching coefficients have widely been explored (Hu and Tian, 2019; Yang et al., 2019; Zhu and Yin, 2009); however, partial observation problems of this type focused on in this paper have not been considered so far.

Based on Wang (2001) and Pham and Tankov (2009) with suitable modifications, we formulate the dynamic programming principle for our problem and derive the optimality equation to find the optimal policy. We show that the optimality equation is a coupled system containing degenerate elliptic partial integro-differential equation and a classical minimization problem. Its solvability is analyzed from a viscosity viewpoint. A difference between our and the previous models is that our model considers a controlled observation process, while the previous ones do not. In addition, we derive an exact smooth solution to a reduced but still nontrivial problem.

We also numerically compute the optimality equation using coefficients estimated from hydrological and biological data at a study site in Japan. The optimality equation is discretized by a high-resolution finite difference scheme based on the Weighted Essentially Non-Oscillatory (WENO) reconstruction (Jiang and Peng, 2000). The WENO construction allows us to mitigate the slow convergence of the conventional monotone schemes (Oberman, 2006). The optimal policy of the algae population dynamics is then computed and it is compared with a similar model based on a less flexible policy. Our contribution thus covers mathematical and numerical analysis of the new stochastic optimal control problem under partial observation and its application.

### 1.4 Organization of this paper

The rest of this paper is as follows. In Section 2, we explain the system dynamics and observation/control processes. The dynamic programming principle and the optimality equation are then presented. The control problem is analyzed mathematically in Section 3. Solvability of the optimality equation is analyzed from a viscosity viewpoint, and a tractable model is discussed. A comparison principle to guarantee uniqueness of the optimality equation is also discussed. In Section 4, the optimality equation is solved numerically. Summary and future perspectives are presented in Section 5. Several technical proofs of the propositions in the main text are summarized in **Appendix A**. The model based on a less flexible management policy is explained in **Appendix B**. The flexible and inflexible models are briefly compared in **Appendix C**.

## 2. Mathematical model

### 2.1 System dynamics and control

We consider spatially-lumped algae population dynamics in a river environment. The time is denoted as $t \geq 0$. The system dynamics are descrined with the three continuous-time stochastic variables: a Markov chain representing river flow regimes, a bounded non-nagative variable representing the population, and a controlled Poisson process representing an observation process. **Figure 1** shows the conceptual diagram of the observation and control dynamics.

#### 2.1.1 Flow regimes

The river flow regimes are exogenously given and are represented by a Markov chain $\alpha = (\alpha_t)_{t \geq 0}$ having $I+1$ regimes $\{R_i\}_{0 \leq i \leq I}$ with an integer $I \geq 0$. For the sake of simplicity, we assume that the Markov chain is time-homogenous. Set $J = \{0,1,2,...,I\}$. For each $i \in J$, the river flow discharge $Q_i > 0$ is assigned to the regime $R_i$. Without any loss of generality, we assume that $\{Q_i\}_{0 \leq i \leq I}$ is strictly increasing. The switching rates from the regimes $R_i$ to $R_j$ is $w_{i,j} \geq 0$. An advantage of using a Markov chain is its simplicity and flexibility to represent a wide variety of flow regimes in a non-parametric manner.

### 2.1.2 Population dynamics

The population dynamics are affected by the river flow regimes and are impulsively controlled by the decision-maker: the environmental manager. In case of the filamentous algae, it is common to identify the population as the algae biomass per unit riverbed area (Yoshioka and Yaegashi, 2018). The population at time $t$ is denoted as $X_t$. Set $X = (X_t)_{t \geq 0}$. We assume that the population is bounded as in the classical logistic models (Lungu and Øksendal, 1997; Grandits et al., 2019; Yoshioka, 2019) and normalize the range of $X$ as $D = [0,1]$.

Without any interventions by the decision-maker, we set the governing SDE of $X$ as

$$\mathrm{d}X_t = \left(g\left(Q_{\alpha_t}, X_t\right) - \eta Q_{\alpha_t}\right) X_t \mathrm{d}t = f\left(Q_{\alpha_t}, X_t\right)\mathrm{d}t,\quad t > 0 \tag{1}$$

with an initial condition $X_{0+} \in D$ ( $X_{s+} = \lim_{t \downarrow 0} X_{s+t}$, $s \geq 0$ ). Here, $\eta > 0$ is a detachment coefficient and $g : \mathbb{R} \times D \to [0,+\infty)$ is globally Lipschitz continuous in $\mathbb{R} \times D$, concave with respect to the second argument, and $g(\cdot,0) = g(\cdot,K) = 0$ with $0 < K \leq 1$ (environmental capacity). By the concavity, we assume $g(\cdot,y) \geq 0$ for $0 \leq y \leq K$ and $g(\cdot,y) < 0$ for $y > K$. We can extend $g$ to a globally Lipschitz continuous function in $\mathbb{R} \times \mathbb{R}$ as $g(\cdot,y) = 0$ for $y < 0$ and $g(\cdot,y) = g(\cdot,1) \leq 0$ for $y > 1$. By the Lipschitz continuity of the coefficients, the SDE (1) is solvable uniquely in a strong sense (Theorem 2.1 of Yin and Zhu, 2009). The uniform boundedness $0 \leq X_t \leq 1$ for $t \geq 0$ a.s. then holds true by a contradiction argument similar to Appendix A of Lungu and Øksendal (1997). We assume that, given $\alpha$, $X_t$ ($t \geq 0$) is increasing with respect to the initial condition $X_{0+}$, which would be a biologically natural requirement.

The SDE (1) is piecewise continuous and subject to the regime-switching coefficients. The coefficient $g\left(Q_{\alpha_t}, X_t\right) X_t$ represents the limited intrinsic growth of the population (Tomlinson et al., 2010; Yoshioka, 2019), while the coefficient $\eta Q_{\alpha_t} X_t$ is the simplest model representing flushing out of the population due to river flows (Kazama and Watanabe, 2018; King et al., 2014; Katz et al., 2018), such that detachment of the algae population due to hydrodynamic force increases as the discharge increases.

***Remark 1***

We can add a diffusion term to the SDE (1). A candidate would be the one proportional to $X_t(1 - X_t)\mathrm{d}B_t$, where $B = (B_t)_{t \geq 0}$ is a 1-D standard Brownian motion. Adding this term would not deteriorate the unique solvability and boundedness of the SDE (1).

### 2.1.3 Observation and intervention

The observation process is discrete in time such that the sequence $\tau = \{\tau_k\}_{k=0,1,2,...}$ of the observation times follows a controlled Poisson process $N = (N_t)_{t \geq 0}$. Without any loss of generality, set $\tau_0 = 0$. We assume that there exist the two observation schemes with a high intensity $\bar{\lambda}$ and a low intensity $\underline{\lambda}$ ($0 < \underline{\lambda} < \bar{\lambda} < +\infty$). We assume that the intensity can be adaptively switched at each $\tau_k$ by the decision-maker. The intensity in the interval $(\tau_k, \tau_{k+1}]$ is denoted as $\lambda^{(k)}$, which equals $\underline{\lambda}$ or $\bar{\lambda}$. Set $\lambda_t = \sum_{k=0}^{\infty} \lambda^{(k)} \chi_{(\tau_k, \tau_{k+1}]}(t)$ for $t \geq 0$, where $\chi_{(\tau_k, \tau_{k+1}]}(t)$ is the indicator function with $\chi_{(\tau_k, \tau_{k+1}]}(t) = 1$ if $t \in (\tau_k, \tau_{k+1}]$ and $\chi_{(\tau_k, \tau_{k+1}]}(t) = 0$ otherwise. Due to the strict positivity of the intensity, we have a.s., $\tau_k < \tau_{k+1}$ ($k = 0,1,2,...$) and $\tau_k \to +\infty$ as $k \to +\infty$.

The decision-maker can impulsively harvest the population at each $\tau_k$ through a cleaning up activity (Ismail and Salim, 2013; Joshi et al., 2012; Tezuka et al., 2014). The harvesting is assumed to be proportional to the population at each $\tau_k$:

$$X_{\tau_k +} = (1 - z_k) X_{\tau_k}, \tag{2}$$

where $z_k \in \{0, \bar{z}\}$ with $0 < \bar{z} < 1$ and $z_0 = 0$ is the harvesting rate at $\tau_k$. The relationship (2) serves as a kind of internal boundary condition at $\tau_k$. In this way, we consistently represent both the observation with ($z_k = \bar{z}$) and without harvesting ($z_k = 0$) at each $\tau_k$.

We assume that $z_k$ is decided based on the information available at $\tau_{k-1}$ or $\tau_k$. We call the former inflexible and the latter flexible, because the decision-making of the latter is based on the latest information and can more flexibly control the dynamics. Both cases would arise depending on the preference of the environmental manager. The main interest of our analysis is the flexible case that fully utilizes the available data at each observation. Later, we compare both cases numerically and compute the value of information obtained through the flexible interventions. The mathematical analysis in the main text focus on the flexible case, and that for the inflexible case is placed in **Appendix B**. They are briefly compared in **Appendix C**.

In practice, harvesting the population can be carried out through cleaning activities of the riverbed. Such activities can be dangerous for the decision-maker (humans) if the river discharge is sufficiently large. In addition, the net population growth can become strictly negative in such flow regimes, implying less necessity of the cleaning up activities. Therefore, we assume that there is an integer $L \geq 0$ such that $z_k = 0$ is preferred if $\alpha_{\tau_k} > L$. This effect is not incorporated into the range of $z_k$, but in the harvesting cost as explained in the next sub-section.

***Remark 2***

There may be no technical difficulty to generalize $\Lambda$ and $Z$ to finite discrete sets. We focus on the binary case to formulate a simpler model.

### 2.1.4 System dynamics

The complete system dynamics including observation and intervention processes are presented. The natural filtration generated by the observation processes up to the time $t \geq 0$ is set as

$$\mathcal{F}_t = \sigma\left\{\left(\tau_j, \alpha_{\tau_j}, X_{\tau_j}\right)_{0 \leq j \leq k}, k = \sup\left\{j : \tau_j \leq t\right\}\right\} \tag{3}$$

because no information arrives except at each observation. Set $\mathcal{F} = (\mathcal{F}_t)_{t \geq 0}$. The filtration is augmented by the null-observable sets. Recall that the sequence $\tau$ of observation times is random but affected by the decisions by the decision-maker. We assume that $\lambda^{(k)}$, which is the intensity during the time interval $(\tau_k, \tau_{k+1}]$, is decided based on the information available at the latest observation (possibly with harvesting) time $t = \tau_k$. Therefore, we assume that $\lambda^{(k)}$ is $\mathcal{F}_{\tau_k +}$-measurable. The harvesting process is also assumed to be decided based on the information available at the latest observation time $t = \tau_k$; namely, we require that $z_k$ is $\mathcal{F}_{\tau_k}$-measurable.

The admissible set of control policies is denoted as $\mathcal{C}$, which is the set of the pair $u = (\lambda, z)$ of the process $\lambda = (\lambda_t)_{t \geq 0}$ and the sequence $z = \{z_k\}_{k=0,1,2,...}$: for each $k = 0,1,2,...$, $\lambda^{(k)}$ and $z_k$ are $\mathcal{F}_{\tau_k +}$-measurable and $\mathcal{F}_{\tau_k}$-measurable, and satisfy $\lambda^{(k)} \in \Lambda = \{\bar{\lambda}, \underline{\lambda}\}$ and $z_k \in Z$. Clearly, $\mathcal{C}$ is non-empty and independent from $(\alpha_{0+}, X_{0+})$.

Finally, the population dynamics given a control in $\mathcal{C}$ are set as

$$\mathrm{d}X_t = f\left(Q_{\alpha_t}, X_t\right)\mathrm{d}t, \quad t \in (\tau_k, \tau_{k+1}), \quad k = 0,1,2,... \tag{4}$$

subject to the initial condition $X_{0+} \in D$ and the impulsive relationship (2) at each $t = \tau_k$. By the definition of the admissible set $\mathcal{C}$ and the population dynamics, it is straightforward to see that the boundedness requirement, namely the property $0 \leq X_t \leq 1$ for $t \geq 0$ a.s., is satisfied without any extra constraints. The initial condition is specified just after the 0th observation time $\tau_0 = 0$ without any loss of generality.

At the end of this sub-section, we show a continuity estimate of the population dynamics. The population $X$ subject to the initial condition $X_{0+} = x \in D$ is denoted as $X^x = \left(X_t^x\right)_{t \ge 0}$. The expectation with respect to $\mathcal{F}$ is denoted as $\mathbb{E}$.

***Proposition 2.1***

*For any $x, y \in D$ and $u \in \mathcal{C}$, we have*

$$\left(\mathbb{E}\left[\left|X_t^x - X_t^y\right|\right]\right)^2 \le \mathbb{E}\left[\left|X_t^x - X_t^y\right|^2\right] \le e^{2\omega t}\left|x - y\right|^2, \quad t \ge 0 \tag{5}$$

*with some $\omega > 0$.*

**(Proof of Proposition 2.1)**

The left inequality is by the classical Jensenõs inequality. Therefore, we focus on the right inequality. Choose $u = (\lambda, z) \in \mathcal{C}$. By (2), at each $\tau_k$:

$$\left|X_{\tau_k+}^x - X_{\tau_k+}^y\right| = (1 - z_k)\left|X_{\tau_k}^x - X_{\tau_k}^y\right| \le \left|X_{\tau_k}^x - X_{\tau_k}^y\right|. \tag{6}$$

This becomes an equality if there is no intervention at $\tau_k$ ($z_k = 0$). Each harvesting is thus non-expansive. Set $\underline{u} = (\lambda, \underline{z}) \in \mathcal{C}$ with $\underline{z}_k = 0$ ($k \ge 0$). Then, we get

$$\left|X_t^x - X_t^y\right|_u \le \left|X_t^x - X_t^y\right|_{\underline{u}}, \quad t \ge 0. \tag{7}$$

This inequality holds true against arbitrary $u \in \mathcal{C}$, indicating that its right-hand side gives an upper threshold value of $\left|X_t^x - X_t^y\right|_u$ for any $u = (\lambda, z) \in \mathcal{C}$. With the control $\underline{u}$, the process $X$ is governed by the SDE (1). By (6), (7), considering the SDE (1) with Lemma 2.14 of Yin and Zhu (2009) obtains the right inequality of (5).

Ĭ

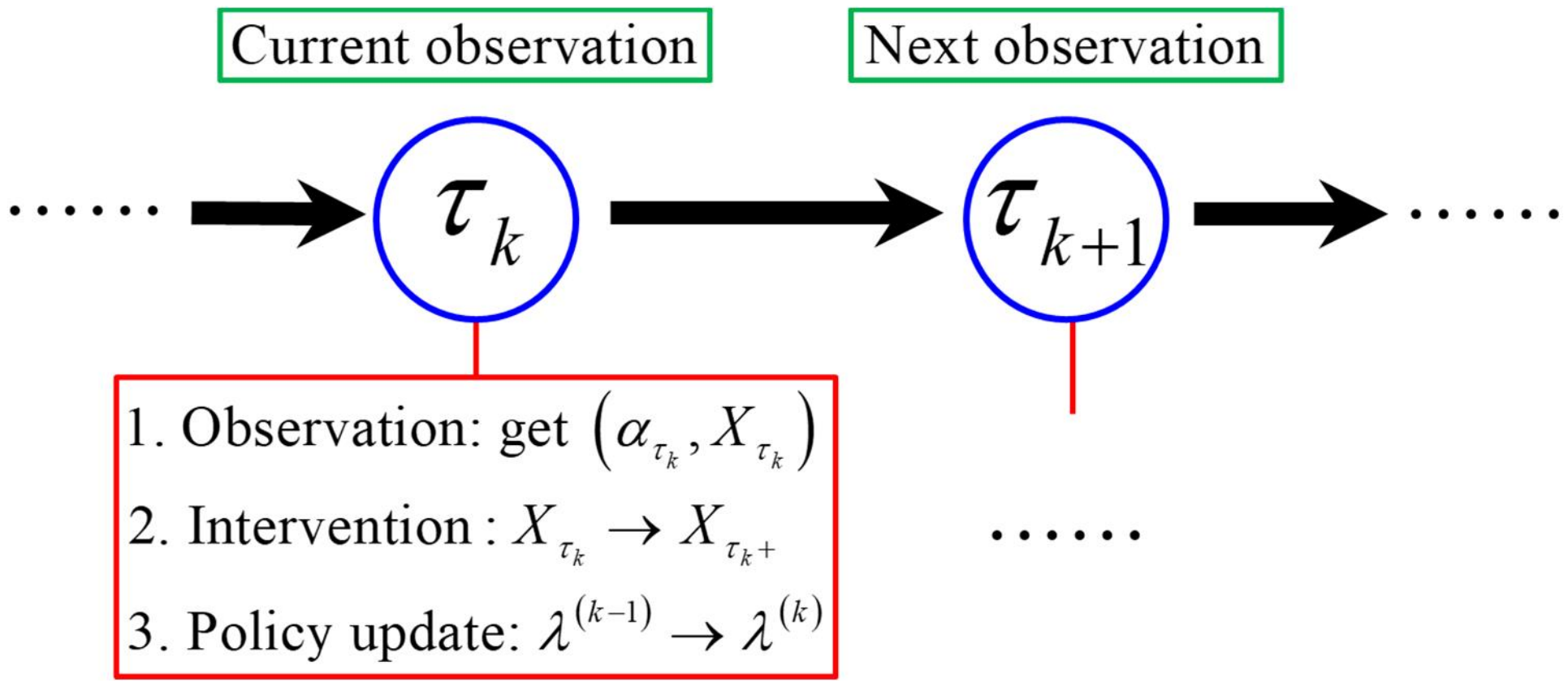

**Figure 1.** Conceptual diagram of the observation and control dynamics.

### 2.2 Performance index

The performance index $\phi: J \times D \times \mathcal{C} \to \mathbb{R}$ to be minimized by choosing the control policy is a conditional expectation containing the cumulative disutility by the population, the observation cost, and the harvesting cost:

$$\phi(i,x;u) = \mathbb{E}^{i,x}\left[\int_0^\infty h(X_s)e^{-\delta s}\mathrm{d}s + \sum_{k=1}^{\infty} e^{-\delta\tau_k}\left(d + K\left(\alpha_{\tau_k}, X_{\tau_k}, z_k\right)\right)\right], \tag{8}$$

where $(\alpha_{0+}, X_{0+}) = (i,x)$ and $\delta > 0$ is the discount rate, $d > 0$ is the observation cost, $h: D \to \mathbb{R}$ is a non-negative increasing Hölder continuous function representing disutility caused by the population, and $K: D \times Z \to \mathbb{R}$ is the harvesting cost:

$$K\left(\alpha_{\tau_k}, X_{\tau_k}, z_k\right) = \begin{cases} 0 & (z = 0) \\ K_0 + K_1\theta\left(X_{\tau_k}\right)z & \left(z = \bar{z},\ \alpha_{\tau_k} \le L\right) \\ P\left(K_0 + K_1\theta\left(X_{\tau_k}\right)z\right) & \left(z = \bar{z},\ \alpha_{\tau_k} > L\right) \end{cases} \tag{9}$$

with $P > 1$ and a function $\theta: D \to [0, +\infty)$ that is Lipschitz continuous in $D$ with $\theta(1) > 0$.

The harvesting cost is not incurred if population is not harvested, while it contains the fixed cost $K_0 \ge 0$ and the proportional cost with $K_1 > 0$ if the population is harvested. No cost is incurred when switching the observation rate $\lambda$. The harvesting cost becomes relatively larger if there is harvesting under a large flow discharge like a flood ($z = \bar{z},\ \alpha_{\tau_k} > L$), because the cleaning activities under such a flow regime can be dangerous and risky. The observation cost $d$ prevents the decision-maker from unconditionally choosing the finer observation scheme.

The value function $\Phi: J \times D \to \mathbb{R}$ is the minimized performance index:

$$\Phi(i,x) = \inf_{u \in \mathcal{C}} \phi(i,x;u). \tag{10}$$

Clearly, $\Phi$ is non-negative in $J \times D$; namely, $\Phi(i,x) \ge 0$ for all $(i,x) \in J \times D$. It is uniformly bounded since

$$\begin{aligned} 0 \le \Phi(i,x) &\le \mathbb{E}^{i,x}\left[\int_0^\infty h(1)e^{-\delta s}\mathrm{d}s + \sum_{k=1}^{\infty} e^{-\delta\tau_k}\left(d + K\left(\alpha_{\tau_k}, 1, 0\right)\right)\right] \\ &\le \frac{h(1)}{\delta} + d\mathbb{E}^{i,x}\left[\sum_{k=1}^{\infty} e^{-\delta\tau_k}\right] \\ &\le \frac{h(1)}{\delta} + d\frac{\bar{\lambda}}{\delta} \end{aligned}. \tag{11}$$

A control $u \in \mathcal{C}$ minimizing $\phi$, assuming the existence, is referred to as an optimal control, and is denoted as $u^*$.

### 2.3 Optimality equation

In this sub-section, firstly, we prove continuity of $\Phi$ for all $\delta > 0$, which will be used in derivation of a dynamic programming principle of our problem. In addition, the continuity result supports that solutions to our problem are (continuous) viscosity solutions.

***Proposition 2.2***

*For each* $i \in J$, $\Phi(i,x)$ *is increasing and continuous with respect to* $x \in D$.

**(Proof of Proposition 2.2)**

Fix some $i \in J$. The increasing nature of $X^x$ with respect to the initial condition $x$ combined with the increasing nature of $h$ and the functional form (9) of $K$ leads to that the $\Phi(i,x)$ is increasing with respect to $x \in D$.

Optimal controls with respect to the initial conditions $(i,x)$ and $(i,y)$ are denoted as $u^x = (\lambda^x, z^x)$ and $u^y = (\lambda^y, z^y)$, respectively. By their optimality, we have

$$\begin{aligned} \Phi(i,x) - \Phi(i,y) &= \inf_{u \in C} \phi(i,x;u) - \inf_{u \in C} \phi(i,y;u) \\ &= \phi(i,x;u^x) - \phi(i,y;u^y) \\ &\le \phi(i,x;u^y) - \phi(i,y;u^y) \\ &= \mathbb{E}\left[ \int_0^\infty \left( h(X_s^x) - h(X_s^y) \right) e^{-\delta s} \mathrm{d}s \right] \\ &\quad + \mathbb{E}\left[ \sum_{k=1}^\infty e^{-\delta \tau_k} \left( K\left(\alpha_{\tau_k}, X_{\tau_k}^x, z_k^y\right) - K\left(\alpha_{\tau_k}, X_{\tau_k}^y, z_k^y\right) \right) \right] \end{aligned} . \tag{12}$$

Here, $u = u^y$ in the last line. The first expectation in the last line of (12) is evaluated with the Hölder continuity of $h$ and the classical Jensenś inequality:

$$\begin{aligned} \mathbb{E}\left[ \int_0^\infty \left( h(X_s^x) - h(X_s^y) \right) e^{-\delta s} \mathrm{d}s \right] &\le \mathbb{E}\left[ \int_0^\infty h_0 \left| X_s^x - X_s^y \right|^\alpha e^{-\delta s} \mathrm{d}s \right] \\ &= h_0 \int_0^\infty \mathbb{E}\left[ \left| X_s^x - X_s^y \right|^\alpha \right] e^{-\delta s} \mathrm{d}s \\ &\le h_0 \int_0^\infty \left( \mathbb{E}\left[ \left| X_s^x - X_s^y \right| \right] \right)^\alpha e^{-\delta s} \mathrm{d}s \end{aligned} \tag{13}$$

with some $h_0 > 0$ and $0 < \alpha \le 1$.

We separately consider the two mutually exclusive cases $\delta - \alpha\omega > 0$ and $\delta - \alpha\omega \le 0$.

***Case 1:*** $\delta - \alpha\omega > 0$

Applying **Proposition 2.1** to (13) yields

$$\int_0^\infty \left(\mathbb{E}\left[\left|X_s^x - X_s^y\right|\right]\right)^\alpha e^{-\delta s}\mathrm{d}s \le \int_0^\infty |x-y|^\alpha e^{-(\delta-\alpha\omega)s}\mathrm{d}s = \frac{1}{\delta-\alpha\omega}|x-y|^\alpha . \tag{14}$$

The first expectation in the last line of (12) is evaluated as

$$\mathbb{E}\left[\sum_{k=1}^{\infty} e^{-\delta\tau_k}\left(K\left(\alpha_{\tau_k}, X_{\tau_k}^x, z_k^y\right) - K\left(\alpha_{\tau_k}, X_{\tau_k}^y, z_k^y\right)\right)\right] \le C_1 \sum_{k=1}^{\infty} \mathbb{E}\left[e^{-\delta\tau_k}\left|X_{\tau_k}^x - X_{\tau_k}^y\right|\right] \tag{15}$$

with a constant $C_1 > 0$. Set $\underline{u} = (\lambda, \underline{z}) \in \mathcal{C}$ with $\underline{z}_k = 0$ ($k \ge 0$). As in the proof of **Proposition 2.1**, we get

$$\mathbb{E}\left[e^{-\delta\tau_k}\left|X_{\tau_k}^x - X_{\tau_k}^y\right|\right] \le \mathbb{E}\left[e^{-\delta\tau_k}\left|X_{\tau_k}^x - X_{\tau_k}^y\right|\right]\Big|_{u=\underline{u}} . \tag{16}$$

Under the control $\underline{u}$, the process $X$ is governed by the SDE (1). By (16), we have

$$\begin{aligned}\mathbb{E}\left[e^{-\delta\tau_k}\left|X_{\tau_k}^x - X_{\tau_k}^y\right|\right] &\le \mathbb{E}\left[e^{-\delta\tau_k}\left|X_{\tau_k}^x - X_{\tau_k}^y\right|\right]_{u=\underline{u}} \\ &= \mathbb{E}\left[e^{-\delta\tau_k}\mathbb{E}\left[\left|X_{\tau_k}^x - X_{\tau_k}^y\right| \middle| \mathcal{F}_{\tau_{k-1}}\right]\right]_{u=\underline{u}} . \\ &\le \mathbb{E}\left[e^{-\delta\tau_k}\mathbb{E}\left[\left|X_{\tau_k}^x - X_{\tau_k}^y\right|^\alpha \middle| \mathcal{F}_{\tau_{k-1}}\right]\right]_{u=\underline{u}}\end{aligned} \tag{17}$$

Here, $x \le x^\alpha$ for any $x \in [0,1]$ and $\alpha \in (0,1]$. Again by **Proposition 2.1** and Jensenós inequality, for $k \ge 1$ we have

$$\begin{aligned}\mathbb{E}\left[e^{-\delta\tau_k}\mathbb{E}\left[\left|X_{\tau_k}^x - X_{\tau_k}^y\right|^\alpha \middle| \mathcal{F}_{\tau_{k-1}}\right]\right]_{u=\underline{u}} &\le \mathbb{E}\left[e^{-\delta\tau_k}\left(\mathbb{E}\left[\left|X_{\tau_k}^x - X_{\tau_k}^y\right| \middle| \mathcal{F}_{\tau_{k-1}}\right]\right)^\alpha\right]_{u=\underline{u}} . \\ &\le \mathbb{E}\left[e^{-\delta\tau_k}e^{\alpha\omega(\tau_k-\tau_{k-1})}\left|X_{\tau_{k-1}}^x - X_{\tau_{k-1}}^y\right|^\alpha\right]_{u=\underline{u}}\end{aligned} \tag{18}$$

The summation in the right-hand side of (15) is denoted as $G$. By (17)-(18), we have

$$\begin{aligned}G &\le \sum_{k=1}^{\infty}\mathbb{E}\left[e^{-\delta\tau_k}e^{\alpha\omega(\tau_k-\tau_{k-1})}\left|X_{\tau_{k-1}}^x - X_{\tau_{k-1}}^y\right|^\alpha\right]\Big|_{u=\underline{u}} \\ &= \mathbb{E}\left[e^{-\delta\tau_1}e^{\alpha\omega(\tau_1-\tau_0)}\left|X_{\tau_0}^x - X_{\tau_0}^y\right|^\alpha\right]\Big|_{u=\underline{u}} + \sum_{k=2}^{\infty}\mathbb{E}\left[e^{-\delta\tau_k}e^{\alpha\omega(\tau_k-\tau_{k-1})}\left|X_{\tau_{k-1}}^x - X_{\tau_{k-1}}^y\right|^\alpha\right]\Big|_{u=\underline{u}} . \\ &= \mathbb{E}\left[e^{-(\delta-\alpha\omega)\tau_1}\right]|x-y|^\alpha + \sum_{k=1}^{\infty}\mathbb{E}\left[e^{-(\delta-\alpha\omega)(\tau_{k+1}-\tau_k)}e^{-\delta\tau_k}\left|X_{\tau_k}^x - X_{\tau_k}^y\right|^\alpha\right]\Big|_{u=\underline{u}}\end{aligned} \tag{19}$$

We evaluate the two terms in the last line of (19). On the first term, we have

$$\mathbb{E}\left[e^{-(\delta-\alpha\omega)\tau_1}\right] \le \frac{\bar{\lambda}}{\bar{\lambda}+\delta-\alpha\omega} < +\infty . \tag{20}$$

The second term is evaluated as

$$
\begin{aligned}
&\sum_{k=1}^{\infty}\mathbb{E}\left[e^{-(\delta-\alpha\omega)(\tau_{k+1}-\tau_k)}e^{-\delta\tau_k}\left|X_{\tau_k}^x-X_{\tau_k}^y\right|^\alpha\right]\Bigg|_{u=\underline{u}}\\
&=\sum_{k=1}^{\infty}\mathbb{E}\left[e^{-\delta\tau_k}\left|X_{\tau_k}^x-X_{\tau_k}^y\right|^\alpha\mathbb{E}\left[e^{-(\delta-\alpha\omega)(\tau_{k+1}-\tau_k)}\middle|\mathcal{F}_{\tau_k}\right]\right]\Bigg|_{u=\underline{u}}\\
&\le\frac{\bar{\lambda}}{\bar{\lambda}+\delta-\alpha\omega}\sum_{k=1}^{\infty}\mathbb{E}\left[e^{-\delta\tau_k}\left|X_{\tau_k}^x-X_{\tau_k}^y\right|^\alpha\right]\Bigg|_{u=\underline{u}}\\
&=\frac{\bar{\lambda}}{\bar{\lambda}+\delta-\alpha\omega}G
\end{aligned}
\tag{21}
$$

since $e^{-(\delta-\alpha\omega)(\tau_{k+1}-\tau_k)}$ is independent from $\mathcal{F}_{\tau_k}$. Substituting (20) and (21) into (17) yields

$$
G\le\frac{\bar{\lambda}}{\bar{\lambda}+\delta-\alpha\omega}|x-y|^\alpha+\frac{\bar{\lambda}}{\bar{\lambda}+\delta-\alpha\omega}G\ \text{ or equivalently }\ G\le\frac{\bar{\lambda}}{\delta-\alpha\omega}|x-y|^\alpha. \tag{22}
$$

Finally, combining (12), (14), (15), and (22) leads to

$$
\Phi(i,x)-\Phi(i,y)\le\frac{h_0+\bar{\lambda}C_1}{\delta-\alpha\omega}|x-y|^\alpha. \tag{23}
$$

Consequently, we get

$$
0\le\Phi(i,x)-\Phi(i,y)\le\frac{h_0+\bar{\lambda}C_1}{\delta-\alpha\omega}|x-y|^\alpha \tag{24}
$$

for all $0\le x,y\le 1$ and completes the proof.

***Case 2:*** $\delta-\alpha\omega\le 0$

We need a stronger boundedness estimate of $\mathbb{E}\left[\left|X_s^x-X_s^y\right|\right]$ to handle the small $\delta$. We already have a.s. $0\le X_s^x,X_s^y\le 1$ and thus $0\le\left|X_s^x-X_s^y\right|\le 1$. Therefore, we get an improved estimate

$$
\mathbb{E}\left[\left|X_s^x-X_s^y\right|\right]\le\min\left\{e^{\omega s}|x-y|,1\right\}. \tag{25}
$$

Following the Proof of Proposition 2.3 of Yoshioka and Tsujimura (2020) for $n=0$, by (25), we get

$$
\mathbb{E}\left[\int_0^\infty\left(h\left(X_s^x\right)-h\left(X_s^y\right)\right)e^{-\delta s}\mathrm{d}s\right]\le h_0F_0\left(|x-y|\right) \tag{26}
$$

with $F_0:D\to\mathbb{R}$ given by

$$
F_0(p)=\begin{cases}\dfrac{1}{\alpha\omega-\delta}\left(\dfrac{\alpha\omega}{\delta}p^{\frac{\delta}{\alpha\omega}}-p\right) & (\delta-\alpha\omega<0)\\[2ex] \dfrac{1}{\delta}p(1-\ln p) & (\delta-\alpha\omega=0)\end{cases}. \tag{27}
$$

This $F_0$ is concave and increasing in $D$ with $F_0(0)=0$ and $F_1(1)<+\infty$.

Secondly, as in Case 1, we choose the no-switching control $u = \underline{u}$. For $k \geq 1$, set

$$G_k = \mathbb{E}\left[ e^{-\delta \tau_k} \left| X^x_{\tau_k} - X^y_{\tau_k} \right| \right] \tag{28}$$

and $G = \sum_{k=1}^{\infty} G_k$. We show that each $G_k$ is uniformly bounded and continuous with respect to $x, y \in D$, with which it follows that $G$ is uniformly continuous with respect to $x, y \in D$. We evaluate the right-hand side of (28) as

$$\begin{aligned}
\mathbb{E}\left[ e^{-\delta \tau_k} \left| X^x_{\tau_k} - X^y_{\tau_k} \right| \right] &= \mathbb{E}\left[ e^{-\delta(\tau_k - \tau_{k-1})} \mathbb{E}\left[ e^{-\delta \tau_{k-1}} \left| X^x_{\tau_k} - X^y_{\tau_k} \right| \middle\| \mathcal{F}_{\tau_{k-1}} \right] \right] \\
&\leq \mathbb{E}\left[ \frac{\bar{\lambda}}{\delta + \bar{\lambda}} \mathbb{E}\left[ e^{-\delta \tau_{k-1}} \left| X^x_{\tau_k} - X^y_{\tau_k} \right| \middle\| \mathcal{F}_{\tau_{k-1}} \right] \right] \\
&\leq \frac{\bar{\lambda}}{\delta + \bar{\lambda}} \mathbb{E}\left[ \mathbb{E}\left[ \left| X^x_{\tau_k} - X^y_{\tau_k} \right| \middle\| \mathcal{F}_{\tau_{k-1}} \right] \right] \\
&= \frac{\bar{\lambda}}{\delta + \bar{\lambda}} \mathbb{E}\left[ \left| X^x_{\tau_k} - X^y_{\tau_k} \right| \right]
\end{aligned} \tag{29}$$

and thus

$$(0 \leq) G_k = \frac{\bar{\lambda}}{\delta + \bar{\lambda}} \mathbb{E}\left[ \left| X^x_{\tau_k} - X^y_{\tau_k} \right| \right] < 1. \tag{30}$$

Therefore, $G_k$ is uniformly bounded. With the help of (25) and the Poisson nature of the observations, for $k \geq 1$, the expectation in (30) is estimated as

$$\begin{aligned}
\mathbb{E}\left[ \left| X^x_{\tau_k} - X^y_{\tau_k} \right| \right] &= C \breve{\mathbb{E}}\left[ \int_{D_k} \lambda^{(0)} \lambda^{(1)} \ldots \lambda^{(k-1)} e_k \left| X^x_{t_k} - X^y_{t_k} \right| \mathrm{d}t_1 \mathrm{d}t_2 \ldots \mathrm{d}t_k \right] \\
&\leq C \breve{\mathbb{E}}\left[ \int_{D_k} \lambda^{(0)} \lambda^{(1)} \ldots \lambda^{(k-1)} e_k \min\left\{ e^{\omega t_k} |x - y|, 1 \right\} \mathrm{d}t_1 \mathrm{d}t_2 \ldots \mathrm{d}t_k \right]
\end{aligned} \tag{31}$$

with a normalization constant $C > 0$, $e_k = \sum_{j=1}^{k} e^{-\lambda^{(j-1)}(\tau_j - \tau_{j-1})}$, $D_k = \left\{ (t_1, t_2, \ldots, t_k) \middle| 0 \leq t_1 \leq t_2, \ldots, \leq t_k \right\}$. Here, $\breve{\mathbb{E}}$ is the expectation with respect to $\alpha$. Each $\lambda^{(j)}$ is random but is either $\bar{\lambda}$ or $\underline{\lambda}$. In particular, there is at least one sequence of $\kappa = \left\{ \kappa_j \right\}_{j=0,1,2,\ldots,k-1}$ each $\kappa_j$ being either $\bar{\lambda}$ or $\underline{\lambda}$ such that

$$\begin{aligned}
&\breve{\mathbb{E}}\left[ \int_{D_k} \lambda^{(0)} \lambda^{(1)} \ldots \lambda^{(k-1)} e_k \min\left\{ e^{\omega t_k} |x - y|, 1 \right\} \mathrm{d}t_1 \mathrm{d}t_2 \ldots \mathrm{d}t_k \right] \\
&\leq \int_{D_k} \kappa_0 \kappa_1 \ldots \kappa_{k-1} \bar{e}_k \min\left\{ e^{\omega t_k} |x - y|, 1 \right\} \mathrm{d}t_1 \mathrm{d}t_2 \ldots \mathrm{d}t_k
\end{aligned} \tag{32}$$

with $\bar{e}_k = \sum_{j=1}^{k} e^{-\kappa_{j-1}(t_j - t_{j-1})}$. There are at most $2^k$ realizations of $\kappa$. From (31) and (32), we get

$$\mathbb{E}\left[\left|X_{\tau_k}^x - X_{\tau_k}^y\right|\right] \le C\int_{D_k} \kappa_0\kappa_1...\kappa_{k-1}\overline{e}_k \min\left\{e^{\omega t_k}\left|x-y\right|,1\right\}\mathrm{d}t_1\mathrm{d}t_2...\mathrm{d}t_k \ . \tag{33}$$

The right-hand side of (33) is an improper integral, which uniformly converges and is continuous with respect to $x, y \in D$ by an application of Theorem 16.10 of Jost (2005), indicating that the last expectation of (31) is also uniformly continuous with respect to $x, y \in D$. In fact, we have

$$C\int_{D_k} \kappa_0\kappa_1...\kappa_{k-1}\overline{e}_k \min\left\{e^{\omega t_k}\left|x-y\right|,1\right\}\mathrm{d}t_1\mathrm{d}t_2...\mathrm{d}t_k \le C\int_{D_k} \kappa_0\kappa_1...\kappa_{k-1}\overline{e}_k \mathrm{d}t_1\mathrm{d}t_2...\mathrm{d}t_k \le 1 \ . \tag{34}$$

By $\left|X_{\tau_k}^x - X_{\tau_k}^y\right| \le 1$, we can directly show the uniform boundedness of $G$ as

$$0 \le G \le \sum_{k=1}^{\infty} G_k \le \sum_{k=1}^{\infty}\left(\frac{\overline{\lambda}}{\delta+\overline{\lambda}}\right)^k < +\infty \ . \tag{35}$$

Then, we can express $G$ as a uniformly continuous function of $0 \le p = \left|x-y\right| \le 1$ as $G = G(p)$ with $G(0) = 0$. In summary, we get the estimate

$$\left|\Phi(i,x) - \Phi(i,y)\right| \le h_0 F_0\left(\left|x-y\right|\right) + C_1 G\left(\left|x-y\right|\right). \tag{36}$$

The right-hand side of (36) is uniformly continuous with respect to $x, y \in D$ and vanishes when $x = y$, showing the continuity of $\Phi$.

□

***Remark 3***

Being different from many of the conventional problems (Øksendal and Sulem, 2019), we do not have to assume a small $\delta$ because of the boundedness of the process $X$. The regularity estimate of $\Phi$ degrades as $\delta$ becomes smaller.

For convenience, for generic $F : J \times D \to \mathbb{R}$, set the intervention operator $\mathcal{M}$ as

$$\mathcal{M}F(i,x) = \inf_{z \in Z}\left\{F\left(i,(1-z)x\right) + d + K(i,x,z)\right\}. \tag{37}$$

Set the auxiliary quantity $\Psi : J \times D \times \Lambda \to \mathbb{R}$:

$$\Psi(i,x,r) = \mathbb{E}^{i,x}\left[\int_0^{\tau_1} h\left(\overline{X}_s\right)e^{-\delta s}\mathrm{d}s + \mathcal{M}\Phi\left(\alpha_{\tau_1}, \overline{X}_{\tau_1}\right)e^{-\delta\tau_1}\right], \tag{38}$$

where $\mathbb{E}^{i,x}$ represents the conditional expectation with respect to $\alpha$ and $N$, under the condition $\left(\alpha_{0+}, \overline{X}_{0+}\right) = (i,x) \in J \times D$. Here, the process $\overline{X} = \left(\overline{X}_t\right)_{t \ge 0}$ is identical to $X$ up to $t = \tau_1$, while no intervention is carried out for $t > \tau_1$.

We state the dynamic programming principle. Its proof is placed in **Appendix A**.

***Proposition 2.3***

*We have the dynamic programming principle*

$$\Phi(i,x) = \inf_{r\in\Lambda} \Psi(i,x,r) = \inf_{r\in\Lambda} \mathbb{E}^{i,x}\left[\int_0^{\tau_1} h(\bar{X}_s)e^{-\delta s}\mathrm{d}s + \mathcal{M}\Phi(\alpha_{\tau_1},\bar{X}_{\tau_1})e^{-\delta\tau_1}\right]. \tag{39}$$

The minimization in (39) is considered over the two cases since $\lambda^{(0)}$ is either $\bar{\lambda}$ or $\underline{\lambda}$. Now, we derive the optimality equation. The right-hand side of (38) is rewritten as

$$\begin{aligned}
&\breve{\mathbb{E}}^{i,x}\left[\int_0^\infty re^{-rt}\left\{\int_0^t h(\bar{X}_s)e^{-\delta s}\mathrm{d}s + \mathcal{M}\Phi(\alpha_t,\bar{X}_t)e^{-\delta t}\right\}\mathrm{d}t\right] \\
&= \breve{\mathbb{E}}^{i,x}\left[\int_0^\infty re^{-rt}\left(\int_0^t h(\bar{X}_s)e^{-\delta s}\mathrm{d}s\right)\mathrm{d}t + \int_0^\infty r\mathcal{M}\Phi(\alpha_t,\bar{X}_t)e^{-(\delta+r)t}\mathrm{d}t\right] \\
&= \breve{\mathbb{E}}^{i,x}\left[\int_0^\infty h(\bar{X}_s)e^{-\delta s}\left(\int_s^\infty re^{-rt}\mathrm{d}t\right)\mathrm{d}s + \int_0^\infty r\mathcal{M}\Phi(\alpha_t,\bar{X}_t)e^{-(\delta+r)t}\mathrm{d}t\right], \\
&= \breve{\mathbb{E}}^{i,x}\left[\int_0^\infty \left\{h(\bar{X}_t) + r\mathcal{M}\Phi(\alpha_t,\bar{X}_t)\right\}e^{-(\delta+r)t}\mathrm{d}t\right]
\end{aligned} \tag{40}$$

where $\breve{\mathbb{E}}^{i,x}$ represents the conditional expectation with respect to $\alpha$. Namely, we get

$$\Psi(i,x,r) = \breve{\mathbb{E}}^{i,x}\left[\int_0^\infty \left\{h(\bar{X}_t) + r\mathcal{M}\Phi(\alpha_t,\bar{X}_t)\right\}e^{-(\delta+r)t}\mathrm{d}t\right]. \tag{41}$$

Given $r\in\Lambda$ and the value function $\Phi$, the right-hand side of (41) is a conditional expectation governed by a degenerate elliptic equation. By the classical Feynman-Kac formula for regime-switching processes (Zhu et al., 2015), we heuristically get the optimality equation as a coupled problem of an integro-differential equation and a static optimization problem:

$$\delta\Psi_i - f(i,x)\frac{\partial\Psi_i}{\partial x} + \sum_{\substack{j=0\\ j\neq i}}^{I} w_{i,j}(\Psi_i - \Psi_j) - h(x) + r(\Psi_i - \mathcal{M}\Phi(i,x)) = 0 \tag{42}$$

and

$$\Phi(i,x) = \inf_{r\in\Lambda}\Psi(i,x,r) \tag{43}$$

with the notation $\Psi_i = \Psi_i(x,r) = \Psi(i,x,r)$. Notice that the system (42)-(43) is not a standard degenerate elliptic form (Øksendal and Sulem, 2019).

From the optimality equation, we can deduce the optimal control at each $\tau_k$ as follows:

$$\lambda^{(k),*} = \arg\min_{r\in\Lambda}\Psi(\alpha_{\tau_k+}, X_{\tau_k+}, r) \text{ and } z_k^* = \arg\min_{z\in Z}\left\{\Phi(\alpha_{\tau_k},(1-z)X_{\tau_k}) + K(\alpha_{\tau_k},X_{\tau_k},z)\right\}. \tag{44}$$

In this sense, we can find the optimal control by solving the optimality equation.

## 3. Mathematical analysis of the optimality equation

In this section, we analyze the optimality equation from the standpoints of an exact smooth solution under a simplified condition and a viscosity solution.

### 3.1 Exact solution

An exact solution to a reduced counterpart of the optimality equation, the system (42) and (43), is presented. We show that the reduced optimality equation still admits a non-trivial structure. The reduced optimality equation is

$$\delta\Psi_i - f_i x\frac{\partial\Psi_i}{\partial x} + \sum_{\substack{j=0\\ j\neq i}}^{I} w_{i,j}\left(\Psi_i - \Psi_j\right) - x + r\left(\Psi_i - \mathcal{M}\Phi(i,x)\right) = 0 \tag{45}$$

for $0 \le i \le I$ ( $I = 1$ ) and $x \ge 0$ with (43), where $f_i$ ( $0 \le i \le I$ ) are constants, $\theta = K = const > 0$, $h = x$. Furthermore, we assume $L = 0$, implying that the $R_1$ is the flood regime in which cleaning up the river is dangerous, while the regime $R_0$ is not.

The reduced counterpart assumes the vanishing observation and fixed costs and has linear drift coefficients and a source term. Due to the linearized nature, the implications obtained from the analysis here would be applicable only to the small population $x$. In addition, we neglect switching of the observation rates and flow regimes. Hence, we consider $r > 0$ as a constant. This means that we can formally write $\Psi_i = \Phi_i$ because the problem no longer has the static optimization part on $r$.

Under the above-presented assumptions, (45) becomes

$$\left(\delta + r + w_{i,1-i}\right)\Phi_i - f_i x\frac{\partial\Phi_i}{\partial x} - w_{i,1-i}\Phi_{1-i} = x + r\mathcal{M}\Phi_i . \tag{46}$$

We guess an exact solution of the form

$$\Phi_i(x) = C_i x \tag{47}$$

with a constant $C_i > 0$ ($i = 0,1$). Substituting (47) into (46) yields

$$\left(\delta - f_i + w_{i,1-i}\right)C_i - w_{i,1-i}C_{1-i} = 1 + r\inf_z\left\{z\left(K^{(i)} - C_i\right)\right\},\ \ i = 0,1 \tag{48}$$

with $K^{(0)} = K$ and $K^{(1)} = PK$. The system (48) can be solved analytically, but our focus is whether we get a solution of the form

$$z^* = \begin{cases}\bar{z} & (i=0)\\ 0 & (i=1)\end{cases}, \tag{49}$$

representing the cleaning up activities only in the regime $R_0$. This can be a rational optimal control such that cleaning up the river environment is avoided if the discharge is large. A straightforward calculation leads to the next proposition.

***Proposition 3.1***

*Set*

$$L = \left(\delta - f_0 + w_{0,1} + r\overline{z}\right)\left(\delta - f_1 + w_{1,0}\right) - w_{0,1} w_{1,0} > 0, \tag{50}$$

$$\overline{C}_0 = \frac{\left(\delta - f_1 + w_{1,0}\right)\left(1 + rK\overline{z}\right) + w_{0,1}}{L}, \text{ and } \overline{C}_1 = \frac{\delta - f_0 + w_{0,1} + r\overline{z} + w_{1,0}\left(1 + rK\overline{z}\right)}{L}. \tag{51}$$

*If*

$$\overline{C}_0 > K \text{ and } \overline{C}_0 \le PK, \tag{52}$$

*then, we get the optimal control* $z^*$ *of the form (49) with* $C_i = \overline{C}_i$ $(i = 0,1)$.

A question is whether we can find parameter values such that the inequalities (52) are satisfied. We can rewrite (52) as the following two inequalities:

$$\left[L - \left(\delta - f_1 + w_{1,0}\right) r\overline{z}\right] K < \delta - f_1 + w_{1,0} + w_{0,1} \tag{53}$$

and

$$P \ge \frac{\delta - f_0 + w_{0,1} + r\overline{z} + w_{1,0}\left(1 + rK\overline{z}\right)}{LK}. \tag{54}$$

The inequality (53) is satisfied if the proportional coefficient $K$ of the harvesting is sufficiently small or the quantity $L - \left(\delta - f_1 + w_{1,0}\right) r\overline{z}$ is not larger than 0. The inequality (54) is then satisfied if we choose a sufficiently large $P$. The exact solution indicates that, under a simplified condition, choosing a larger $P$ indeed leads to an optimal control of the form (49).

For more complicated cases, it seems to be hopeless to get an exact solution, motivating us to employing a numerical method. In Section 4, we numerically examine whether the similar optimal controls exist in the original problem.

### 3.2 Viscosity solutions approach

We show that solutions to the optimality equation are characterized in a continuous viscosity sense. In **Definition 3.1**, ñcontinuous functionsò means functions continuous with respect to the second argument $x \in D$. Similarly, ñcontinuously differentiable functionsò means functions continuously differentiable with respect to the second argument $x \in D$. For definitions of viscosity solutions, see the literature (Crandall et al., 1992; Gassiat et al., 2014; Pemy et al., 2008).

***Definition 3.1***

*A pair* $(\Psi, \Phi)$ *of continuous functions* $\Phi : J \times D \to \mathbb{R}$ *and* $\Psi : \Omega \to \mathbb{R}$ *is a viscosity solution* *to the optimality equation, the system* (42) *and* (43)*, if:*

*(a) Sub-solution:* $\Phi(i,x) \le \inf_r \Psi(i,x,r)$ *for all* $(i,x) \in J \times D$*, and for each* $(i,r) \in J \times \Lambda$ :

$$\begin{aligned} &\delta\Psi(i,\bar{x},r) - f(i,\bar{x})\frac{\partial \varphi}{\partial x}(\bar{x}) + \sum_{\substack{j=0 \\ j\neq i}}^{I} w_{i,j}\left(\Psi(i,\bar{x},r) - \Psi(j,\bar{x},r)\right) - h(\bar{x}) \\ &+ r\left[\Psi(i,\bar{x},r) - \mathcal{M}\Phi(i,\bar{x})\right] \le 0 \end{aligned} \tag{55}$$

*for any continuously differentiable* $\varphi : D \to \mathbb{R}$ *and* $(\bar{i},\bar{x}) \in J \times D$ *that is a local maximum of* $\Psi(\cdot,\cdot,r) - \varphi$.

*(b) Super-solution:* $\Phi(i,x) \ge \inf_r \Psi(i,x,r)$ *for all* $(i,x) \in J \times D$*, and for each* $(i,r) \in J \times \Lambda$ :

$$\begin{aligned} &\delta\Psi(i,\bar{x},r) - f(i,\bar{x})\frac{\partial \varphi}{\partial x}(\bar{x}) + \sum_{\substack{j=0 \\ j\neq i}}^{I} w_{i,j}\left(\Psi(i,\bar{x},r) - \Psi(j,\bar{x},r)\right) - h(\bar{x}) \\ &+ r\left[\Psi(i,\bar{x},r) - \mathcal{M}\Phi(i,\bar{x})\right] \ge 0 \end{aligned} \tag{56}$$

*for any continuously differentiable* $\varphi : D \to \mathbb{R}$ *and* $(\bar{i},\bar{x}) \in J \times D$ *that is a local minimum of* $\Psi(\cdot,\cdot,r) - \varphi$.

We show that the pair $(\Psi,\Phi)$ is the unique viscosity solution. To show this, we firstly show that $(\Psi,\Phi)$ is a viscosity solution, and secondly show that the optimality equation admits at most one viscosity solution.

***Proposition 3.2***

*The pair* $(\Psi,\Phi)$ *is a viscosity solution to the optimality equation.*

**(Proof of Proposition 3.2)**

Firstly, $\Phi$ is continuous in $D$. In addition, $\Phi(i,x) = \inf_r \Psi(i,x,r)$ for all $(i,x) \in J \times D$. Given the continuity of $\Phi$, $\Psi$ is a conditional expectation for each $r \in \Lambda$ with respect to the state variables $(\alpha, \bar{X})$. The standard dynamic programming principle in this context yields the viscosity property of $\Psi$ and hence that of the pair $(\Psi,\Phi)$.

Ĭ

***Proposition 3.3***

*For any viscosity sub-solution* $(\underline{\Psi},\underline{\Phi})$ *and viscosity super-solution* $(\bar{\Psi},\bar{\Phi})$*, we have* $\bar{\Psi} \ge \underline{\Psi}$ *in* $\Omega$ *and* $\bar{\Phi} \geq \underline{\Phi}$ *in* $D$.

**(Proof of Proposition 3.3)**

We apply the classical doubling of variables technique (Crandall et al., 1992). A non-standard part is here the treatment of the non-local terms $\mathcal{M}\Phi_i$. We show that these terms can be handled by considering the definition of the operator $\mathcal{M}$. Notice that the domain $\Omega$ is compact, and the optimality equation is satisfied at each point in the domain without resorting to any state-constraints like Katsoulakis (1994), Pham and Tankov (2009), and Yoshioka and Yoshioka (2020).

Assume that there is a constant $M > 0$ such that

$$M = \sup_{i,x,r} \left(\underline{\Psi} - \bar{\Psi}\right)(i,x,r) = \left(\underline{\Psi} - \bar{\Psi}\right)\left(\breve{i},\breve{x},\breve{r}\right) > 0 . \tag{57}$$

We show that this is impossible. Set the auxiliary function

$$\theta(i,x,y,r) = \underline{\Psi}(i,x,r) - \bar{\Psi}(i,y,r) - \frac{(x-y)^2}{2\varepsilon} . \tag{58}$$

with $\varepsilon > 0$. Based on the standard doubling the variables technique, it follows that the maximum of $\theta$ is attained at some $(i_\varepsilon, x_\varepsilon, y_\varepsilon, r_\varepsilon) \in J \times D \times D \times \Lambda$. The maximizer admits a sub-sequence, with an abuse of notations, such that

$$i_\varepsilon \to \breve{i},\ \ r_\varepsilon \to \breve{r},\ \ (x_\varepsilon, y_\varepsilon) \to (\breve{x},\breve{x}) \text{, and } \frac{(x_\varepsilon - y_\varepsilon)^2}{2\varepsilon} \to 0 \ \text{ as } \ \varepsilon \to +0 . \tag{59}$$

Without any loss of generality, it is sufficient to consider a small $\varepsilon > 0$ such that $i_\varepsilon = \breve{i}$ and $r_\varepsilon = \breve{r}$ because $J$ and $\Lambda$ are finite.

Set $p_\varepsilon = \dfrac{x_\varepsilon - y_\varepsilon}{\varepsilon}$. We see that $\underline{\Psi}\left(\breve{i},x,\breve{r}\right) - \bar{\Psi}\left(\breve{i},y_\varepsilon,\breve{r}\right) - \dfrac{(x-y_\varepsilon)^2}{2\varepsilon}$ is maximized at $x = x_\varepsilon$ and $-\left[\underline{\Psi}\left(\breve{i},x_\varepsilon,\breve{r}\right) - \bar{\Psi}\left(\breve{i},y,\breve{r}\right) - \dfrac{(x_\varepsilon - y)^2}{2\varepsilon}\right]$ is minimized at $y = y_\varepsilon$. The functions $\underline{\varphi}(i,x) = \bar{\Psi}\left(\breve{i},y_\varepsilon,\breve{r}\right) + \dfrac{(x-y_\varepsilon)^2}{2\varepsilon}$ and $\bar{\varphi}(i,x) = \underline{\Psi}\left(\breve{i},x_\varepsilon,\breve{r}\right) - \dfrac{(x_\varepsilon - y)^2}{2\varepsilon}$ can therefore be used as test functions for the sub- and super-solutions. Then, we get the inequalities

$$\begin{aligned}
&\delta\underline{\Psi}\left(\breve{i},x_\varepsilon,\breve{r}\right) - f\left(\breve{i},x_\varepsilon\right)p_\varepsilon + \sum_{\substack{j=0\\ j\neq \breve{i}}}^{I} w_{\breve{i},j}\left(\underline{\Psi}\left(\breve{i},x_\varepsilon,\breve{r}\right) - \underline{\Psi}\left(j,x_\varepsilon,\breve{r}\right)\right) - h(x_\varepsilon) \\
&+\breve{r}\left[\underline{\Psi}\left(\breve{i},x_\varepsilon,\breve{r}\right) - \mathcal{M}\underline{\Phi}\left(\breve{i},x_\varepsilon\right)\right] \leq 0
\end{aligned} \tag{60}$$

and

$$\delta\overline{\Psi}\left(\breve{E},y_\varepsilon,\breve{E}\right)-f\left(\breve{E},x_\varepsilon\right)p_\varepsilon+\sum_{\substack{j=0\\ j\neq \breve{E}}}^{I}w_{\breve{E},j}\left(\overline{\Psi}\left(\breve{E},y_\varepsilon,\breve{E}\right)-\overline{\Psi}\left(j,y_\varepsilon,\breve{E}\right)\right)-h\left(y_\varepsilon\right)$$
$$+\breve{E}\left[\overline{\Psi}\left(\breve{E},y_\varepsilon,\breve{E}\right)-\mathcal{M}\overline{\Phi}\left(\breve{E},y_\varepsilon\right)\right]\geq 0 \quad . \tag{61}$$

Combining (60) and (61) yields

$$\delta\underline{\Psi}\left(\breve{E},x_\varepsilon,\breve{E}\right)-f\left(\breve{E},x_\varepsilon\right)p_\varepsilon+\sum_{\substack{j=0\\ j\neq \breve{E}}}^{I}w_{\breve{E},j}\left(\underline{\Psi}\left(\breve{E},x_\varepsilon,\breve{E}\right)-\underline{\Psi}\left(j,x_\varepsilon,\breve{E}\right)\right)-h\left(x_\varepsilon\right)$$
$$+\breve{E}\left[\underline{\Psi}\left(\breve{E},x_\varepsilon,\breve{E}\right)-\mathcal{M}\underline{\Phi}\left(\breve{E},x_\varepsilon\right)\right] \quad . \tag{62}$$
$$-\begin{pmatrix}\delta\overline{\Psi}\left(\breve{E},y_\varepsilon,\breve{E}\right)-f\left(\breve{E},x_\varepsilon\right)p_\varepsilon+\sum_{\substack{j=0\\ j\neq \breve{E}}}^{I}w_{\breve{E},j}\left(\overline{\Psi}\left(\breve{E},y_\varepsilon,\breve{E}\right)-\overline{\Psi}\left(j,y_\varepsilon,\breve{E}\right)\right)-h\left(y_\varepsilon\right)\\ +\breve{E}\left[\overline{\Psi}\left(\breve{E},y_\varepsilon,\breve{E}\right)-\mathcal{M}\overline{\Phi}\left(\breve{E},y_\varepsilon\right)\right]\end{pmatrix}\leq 0$$

Rearranging the equation and taking the limit $\varepsilon\to+0$ with the help of (59) and the continuity of $h$ yields

$$\delta\left(\underline{\Psi}\left(\breve{E},\breve{E},\breve{E}\right)-\overline{\Psi}\left(\breve{E},\breve{E},\breve{E}\right)\right)$$
$$\leq-\sum_{\substack{j=0\\ j\neq \breve{E}}}^{I}w_{\breve{E},j}\left(\underline{\Psi}\left(\breve{E},\breve{E},\breve{E}\right)-\underline{\Psi}\left(j,x,\breve{E}\right)\right)+\sum_{\substack{j=0\\ j\neq \breve{E}}}^{I}w_{\breve{E},j}\left(\overline{\Psi}\left(\breve{E},\breve{E},\breve{E}\right)-\overline{\Psi}\left(j,\breve{E},\breve{E}\right)\right) \tag{63}$$
$$-\breve{E}\left(\underline{\Psi}\left(\breve{E},\breve{E},\breve{E}\right)-\overline{\Psi}\left(\breve{E},\breve{E},\breve{E}\right)\right)+\breve{E}\left(\mathcal{M}\underline{\Phi}\left(\breve{E},\breve{E}\right)-\mathcal{M}\overline{\Phi}\left(\breve{E},\breve{E}\right)\right)$$

or equivalently

$$\left(\frac{\delta}{\breve{E}}+1\right)M\leq-\frac{1}{\breve{E}}\sum_{\substack{j=0\\ j\neq \breve{E}}}^{I}w_{\breve{E},j}\left[\underline{\Psi}\left(\breve{E},\breve{E},\breve{E}\right)-\overline{\Psi}\left(\breve{E},\breve{E},\breve{E}\right)-\left(\underline{\Psi}\left(j,x,\breve{E}\right)-\overline{\Psi}\left(j,\breve{E},\breve{E}\right)\right)\right]$$
$$+\mathcal{M}\underline{\Phi}\left(\breve{E},\breve{E}\right)-\mathcal{M}\overline{\Phi}\left(\breve{E},\breve{E}\right) \quad . \tag{64}$$

By the maximization property of $\left(\breve{E},\breve{E},\breve{E}\right)$, we have

$$\underline{\Psi}\left(\breve{E},\breve{E},\breve{E}\right)-\overline{\Psi}\left(\breve{E},\breve{E},\breve{E}\right)-\left(\underline{\Psi}\left(j,x,\breve{E}\right)-\overline{\Psi}\left(j,\breve{E},\breve{E}\right)\right)\geq 0 \tag{65}$$

for all $j\in I$, meaning that the first term in the right-hand side of (64) is not positive. There exist $\overline{z},\underline{z}\in Z$ such that

$$\mathcal{M}\underline{\Phi}\left(\breve{E},\breve{E}\right)=\underline{\Phi}\left(\breve{E},\breve{E}\left(1-\underline{z}\right)\right)+K\left(\breve{E},\breve{E},\underline{z}\right),\quad \mathcal{M}\overline{\Phi}\left(\breve{E},\breve{E}\right)=\overline{\Phi}\left(\breve{E},\breve{E}\left(1-\overline{z}\right)\right)+K\left(\breve{E},\breve{E},\overline{z}\right). \tag{66}$$

We get

$$
\begin{aligned}
\mathcal{M}\underline{\Phi}\left(\check{E},\check{E}\right)-\mathcal{M}\overline{\Phi}\left(\check{E},\check{E}\right) &= \underline{\Phi}\left(\check{E},\check{E}(1-\underline{z})\right)+K\left(\check{E},\check{E},\underline{z}\right)-\left[\overline{\Phi}\left(\check{E},\check{E}(1-\overline{z})\right)+K\left(\check{E},\check{E},\overline{z}\right)\right] \\
&\le \underline{\Phi}\left(\check{E},\check{E}(1-\overline{z})\right)+K\left(\check{E},\check{E},\overline{z}\right)-\left[\overline{\Phi}\left(\check{E},\check{E}(1-\overline{z})\right)+K\left(\check{E},\check{E},\overline{z}\right)\right]. \\
&= \underline{\Phi}\left(\check{E},\check{E}(1-\overline{z})\right)-\overline{\Phi}\left(\check{E},\check{E}(1-\overline{z})\right)
\end{aligned}
\tag{67}
$$

By the definitions of the sub- and super-solutions, we can proceed as

$$
\begin{aligned}
\mathcal{M}\underline{\Phi}\left(\check{E},\check{E}\right)-\mathcal{M}\overline{\Phi}\left(\check{E},\check{E}\right) &\le \underline{\Phi}\left(\check{E},\check{E}(1-\overline{z})\right)-\overline{\Phi}\left(\check{E},\check{E}(1-\overline{z})\right) \\
&\le \inf_{r\in\Lambda}\underline{\Psi}\left(\check{E},\check{E}(1-\overline{z}),r\right)-\overline{\Phi}\left(\check{E},\check{E}(1-\overline{z})\right) \\
&\le \inf_{r\in\Lambda}\underline{\Psi}\left(\check{E},\check{E}(1-\overline{z}),r\right)-\inf_{r\in\Lambda}\overline{\Psi}\left(\check{E},\check{E}(1-\overline{z}),r\right) \\
&= \underline{\Psi}\left(\check{E},\check{E}(1-\overline{z}),\underline{r}\right)-\overline{\Psi}\left(\check{E},\check{E}(1-\overline{z}),\overline{r}\right)
\end{aligned}
\tag{68}
$$

with some $\overline{r},\underline{r}\in\Lambda$ because $\Lambda$ is finite. By the minimizing property of $\underline{r}$, we get

$$
\begin{aligned}
\underline{\Psi}\left(\check{E},\check{E}(1-\overline{z}),\underline{r}\right)-\overline{\Psi}\left(\check{E},\check{E}(1-\overline{z}),\overline{r}\right) &\le \underline{\Psi}\left(\check{E},\check{E}(1-\overline{z}),\overline{r}\right)-\overline{\Psi}\left(\check{E},\check{E}(1-\overline{z}),\overline{r}\right) \\
&= \left(\underline{\Psi}-\overline{\Psi}\right)\left(\check{E},\check{E}(1-\overline{z}),\overline{r}\right) \\
&\le M
\end{aligned}
\tag{69}
$$

Consequently, we arrive at the inequality

$$
\left(\frac{\delta}{\check{E}}+1\right)M\le M\ ,
\tag{70}
$$

leading to the contradiction $M\le 0$.

Finally, by the definition we have $\overline{\Phi}\ge\inf_r\overline{\Psi}$ and $\underline{\Phi}\le\inf_r\underline{\Psi}$ in $D$. We get the desired result because $\overline{\Phi}\ge\inf_r\overline{\Psi}\ge\inf_r\underline{\Psi}\ge\underline{\Phi}$ and any viscosity solution is a viscosity sub-solution as well as a viscosity super-solution.

Ĭ

***Remark 4***

Following Calder (2018), the proof can be extended to cases where $h$ is piecewise smooth having a finite number of discontinuous points inside $h$. Note that existence of viscosity solutions should be separately analyzed in such cases. A problem having a degenerate diffusion coefficient can be analogously handled with the Ishiiȏs lemma (Crandall et al., 1992).

## 4. Numerical computation

### 4.1 Finite difference scheme

The local Lax-Friedrichs (LLxF) scheme (Jiang and Peng, 2000) combined with the recursion scheme of Pham and Tankov (2008) has been utilized to numerical computation of the optimality equation. The LLxF scheme is one of the simplest monotone finite differences scheme that can compute uniformly bounded, non-oscillatory, and thus stable numerical solutions under some stability condition. Consistency of the scheme is valid for sufficiently smooth solutions.

In our numerical computation, the following degenerate parabolic counterpart is computed, and the value function and optimal controls of the original infinite-horizon problem are approximately obtained at some sufficiently large terminal time $t = T >> 1$:

$$\frac{\partial \Psi_i}{\partial t} + \delta \Psi_i - f(i,x)\frac{\partial \Psi_i}{\partial x} + \sum_{\substack{j=0 \\ j \neq i}}^{I} w_{i,j}\left(\Psi_i - \Psi_j\right) - h(x) + r\left(\Psi_i - \mathcal{M}\Phi_i\right) = 0, \tag{71}$$

$$\Phi(t,i,x) = \inf_{r \in \Lambda} \Psi(t,i,x,r) \tag{72}$$

subject to an initial condition $\Psi_i = 0$ for $t = 0$. A similar discretization applies to the inflexible case as well. The forward Euler time discretization with a sufficiently small time step size $\Delta t > 0$ is used in the temporal integration, with which we can approximate the viscosity solution (Oberman, 2006).

The LLxF scheme is provably convergent, while its computational accuracy may not be sufficient because it has a low accuracy such that sharp profiles of solutions to degenerate elliptic, parabolic, and hyperbolic problems may be diffused out. Our problem may not be an exception as the sharp solution profiles in the computational results demonstrated later suggest. Therefore, we consider its third-order WENO counterpart (Jiang and Peng, 2000). A key ingredient of the WENO is reconstructing numerical solutions in each cell, so that their sharp transitions are accurately reproduced. A drawback is that the resulting scheme is not necessarily monotone, implying that numerical solutions may not converge toward the appropriate viscosity solutions. Nevertheless, it has been empirically found that the WENO can accurately compute viscosity solutions to degenerate parabolic and hyperbolic equations (Jiang and Peng, 2000; Yoshioka et al., 2020; Yoshioka and Yoshioka, 2020).

***Remark 5***

If one considers population dynamics having a diffusion term like that in **Remark 1**, then he/she can use a semi-Lagrangian scheme (Picarelli and Reisinger, 2020).

### 4.2 Demonstrative application

#### 4.2.1 Parameter identification

The hydrological and biological model parameter has been identified based on the operation data of Obara Dam in Hii River, Shimane Prefecture, Japan. The dam-downstream reach has been suffering from the algae bloom of Periphyton after the construction of the dam. For detailed hydrological information of the river and its watershed, see Ide et al. (2019), Yoshioka and Yaegashi (2018), Tabayashi et al. (2017). The local fisheries cooperative has concerns on bloom of green filamentous algae in the downstream river of the dam, but is facing with gradual decrease of the total number of members and their aging, leading to serious lack of the workforce to monitor the river environment. This problem background motivates us to consider cost-effective river environmental management under the partial observations.

An hourly discharge data from the dam is available from Apr 2016 except for several missing data points. In this section, the empirical Markov chain $\alpha$ is identified for the calendar years 2017, 2018, and 2019, assuming that the Markov chain is homogeneous in time, for the sake of analytical simplicity. This seems to be a critical assumption, but can be justified if the discount rate $\delta$ is specified to be not smaller than $O\left(10^{-1}\right)$ (1/day) so that the decision-making based on the current information can be effectively based on the time scale of at most $\delta^{-1}$ (day): at most a seasonal scale.

In the study area, it has been found that the population is not observed for a large discharge such that $O\left(10^{1}\right)$. Based on this finding, we set $I = 40$ and the representative discharge $Q_i = 0.5 + 1.25i$ (m$^3$/s) for $R_i$ ( $i = 0,1,2,...,40$ ) because the regimes having a relatively large discharge is not of high importance when considering the population dynamics.

**Figures 2** through **4** show the estimated Markov chains in the years 2017, 2018, 2019, respectively. We can see that the Markov chain in 2019 is clearly different from that in 2017 and 2018. In fact, the entropy of the Markov chain based on the natural logarithm (Eq. (5.2.8) of Stratonovich (2020)) to measure randomness of the flow regimes in each year is 0.2605 (2017), 0.2606 (2018), and 0.08036 (2019), respectively. Therefore, the Markov chain for 2019 is indeed different from those of 2017 and 2018 from the entropic viewpoint. Considering the estimation results, we only use the data in 2018 and 2019. In summary, the estimated results clearly imply that the discharges observed in 2019 is smaller than the others.

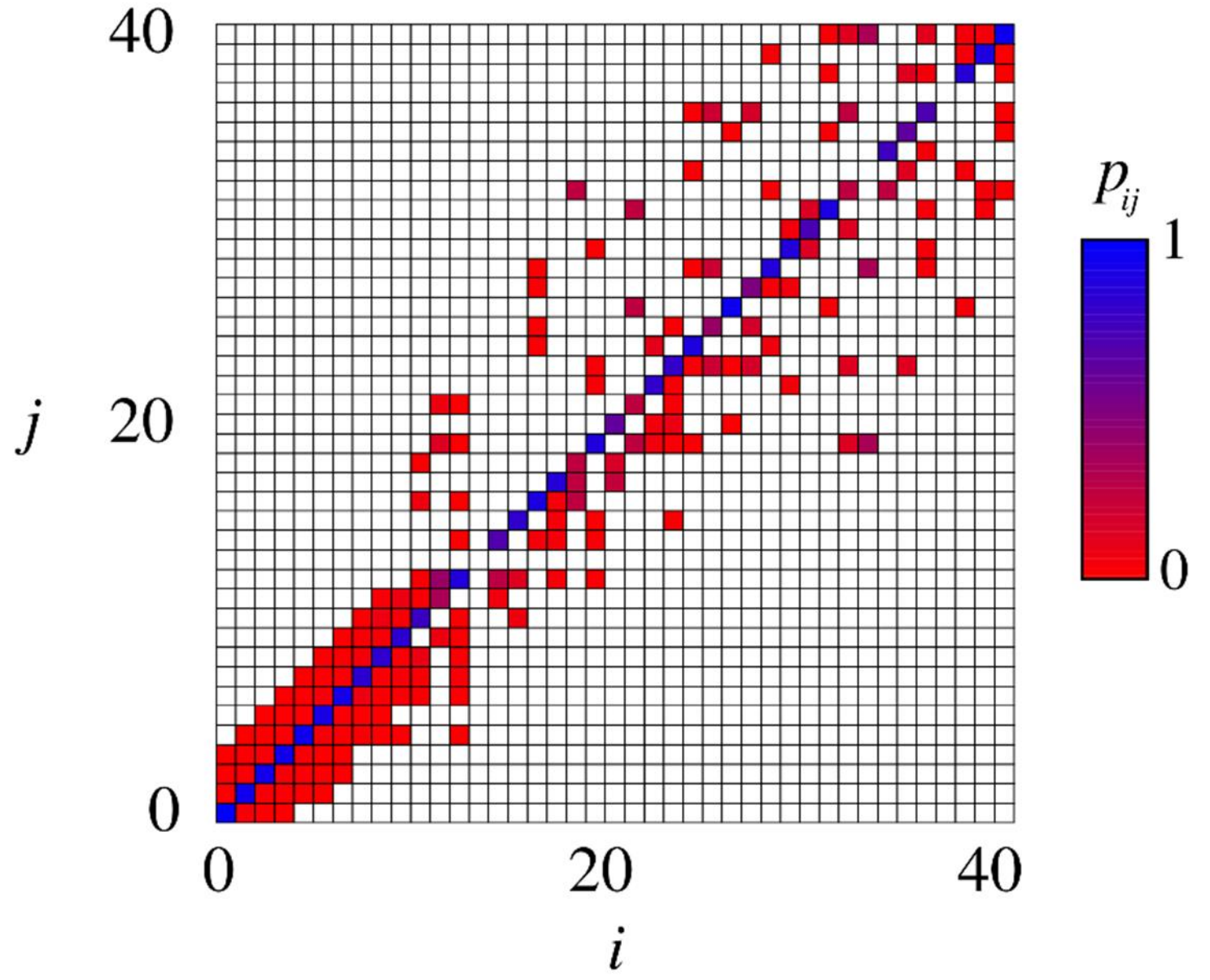

**Figure 2.** Markov chain of the discharge time-series in 2017. The white color represents $p_{ij} = 0$.

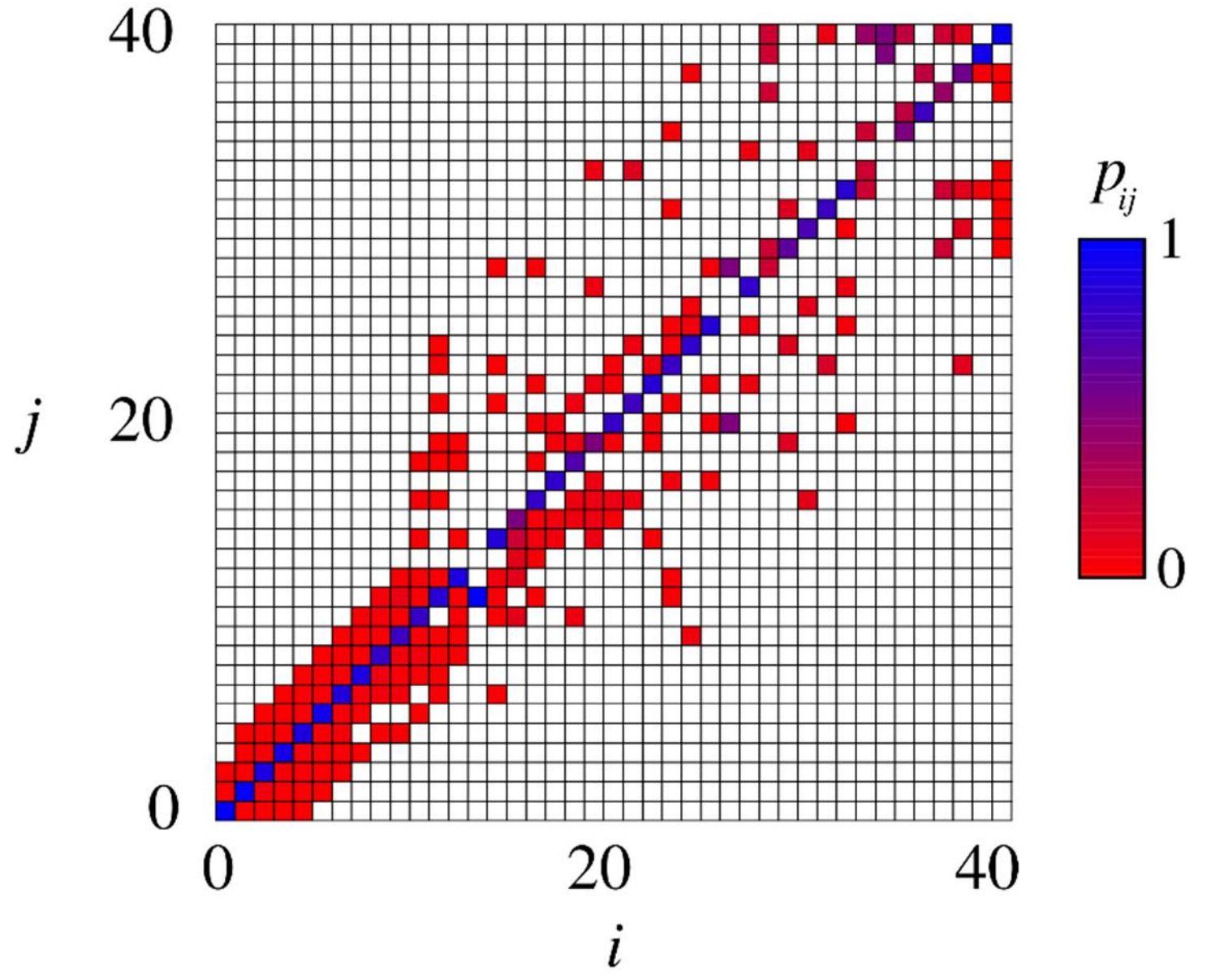

**Figure 3.** Markov chain of the discharge time-series in t 2018. The white color represents $p_{ij} = 0$.

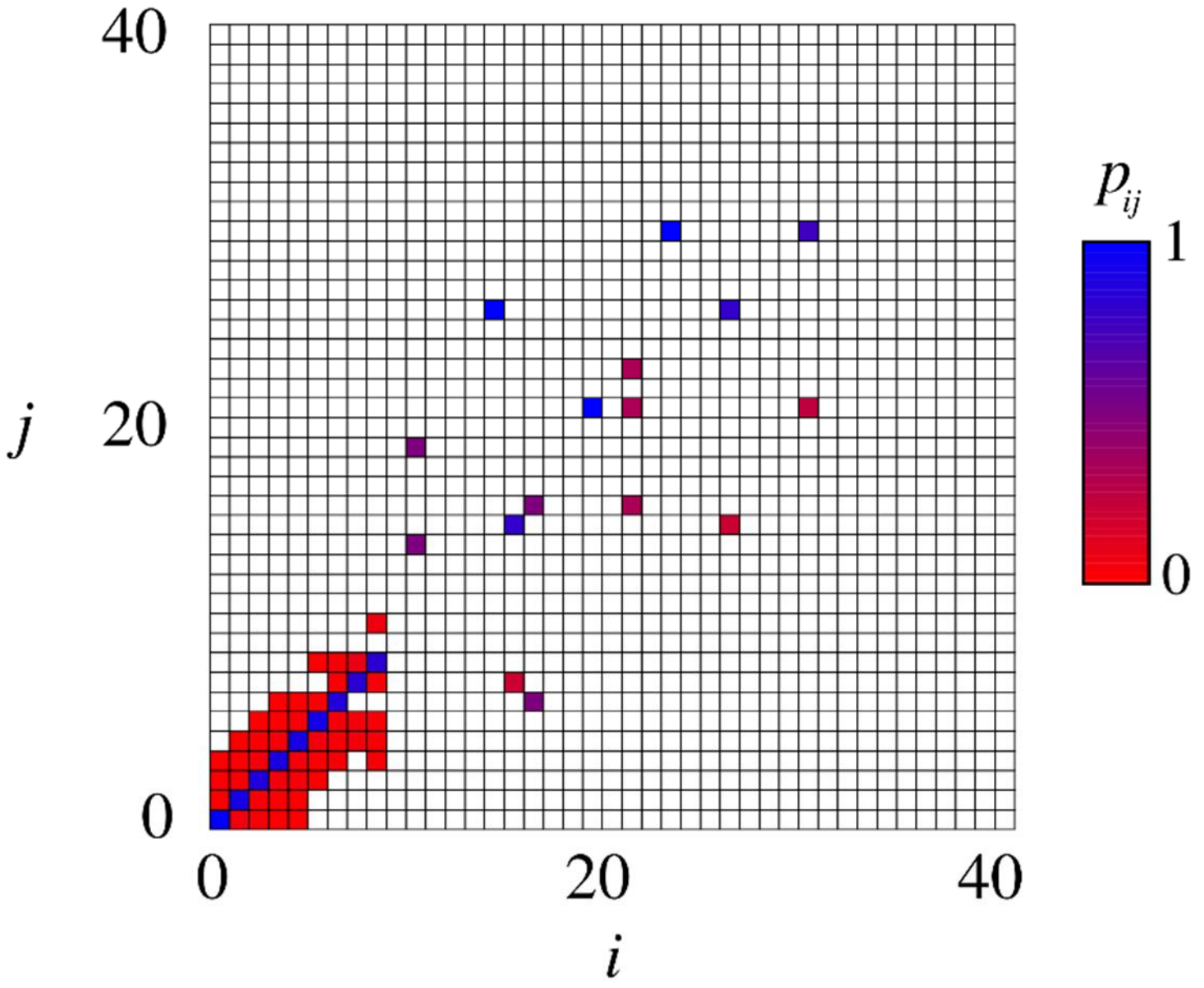

**Figure 4.** Markov chain of the discharge time-series in 2019. The white color represents $p_{ij} = 0$.

Model parameters of the population dynamics are specified as well. Based on Yoshioka (2019), the growth rate $f$ is set as the discharge-dependent logistic model with a normalization

$$f(i,x) = \mu x\left(1 - \frac{x}{c_i}\right) - \eta Q_i x \quad \text{with} \quad 0 < c_i = \frac{aQ_i + b}{aQ_J + b} \leq 1 \tag{73}$$

and constants $\mu > 0$, $a \geq 0$, $b > 0$. This is the simplest logistic model consistent with the concave and unimodal nature of the discharge-dependent stationary population of benthic algae and periphyton, if the regime switching rates $w_{i,j}$ are not so large (Yoshioka, 2019). This model assumes that increasing the river discharge not only leads to larger population decay, but also a larger environmental capacity. The former is due to increasing the hydrodynamic force and physical disturbance acting on the population, while the latter is due to increasing the nutrient transport from the upstream. The case $a = 0$ corresponds to the classical logistic model with the constant capacity $c_i = 1$ and the linear decay $-\eta Q_i x$. If $a > 0$ and $\eta Q_i < \mu$, then solving $f(i,x) = 0$ for $x$ yields $x \propto (\mu - \eta Q_i)(aQ_i + b)$; a concave and unimodal function of the discharge. The model biological parameters have been set as $a = 0.2/Q_I$ (m$^3$/s) and $b = 0.8$ so that $aQ_I + b = 1$, $\mu = 0.5$ (1/day), and $\eta = 0.07$.

We assume that there exist the two options for the observation schemes: the dense scheme with $\lambda = 1/3$ (1/day) and the coarse with $\lambda = 1/10$ (1/day). We set $\delta = 0.2$ (1/day). The harvesting cost is set as $K = K_0 + K_1 xz$ if $z = \bar{z}$ and $K = 0$ otherwise, by which the fixed cost is $K_0$ incurred and the cost $K_1 xz$ with $K_1 > 0$ proportional to the harvested population is incurred at each harvesting. If $i > L = 16$, then $K_0, K_1$ is replaced by $PK_0, PK_1$ with $P = 50$ which we preliminary found that not harvesting the population for $i > L$ is optimal: $z^*(i,x) = 0$ for $i > L$. In our computation, we set $K_0 = 0.15$ and $K_1 = 0.25$. We set $\bar{z} = 0.50$. In addition, set the observation cost $d = 0.1$, implying that it is smaller than the fixed harvesting cost. The disutility $h$ is parameterized as $h = x^m$ with $m > 0$. We set $m = 2$, meaning that the disutility is a convex function of the population.

The finite difference scheme is implemented using a uniform computational grid in both space and time. The total number of vertices used in the space is 401, meaning the spatial increment of 0.0025. The temporal increment for the temporal integration is set as a sufficiently small time step 0.0003 (day) with the sufficiently large terminal time 365/4 (day). The specified parameter values in this sub-section are used unless otherwise specified.

#### 4.2.2 Computational results

**Comparison between 2018 and 2019**

The optimal policies between the different years are numerically analyzed in this section. **Figures 5** and **6** show the value functions and optimal controls ( $\Phi, z^{*}, \lambda^{*}$ ) in 2018 and 2019. The computational results show that the value functions are successfully computed without any spurious oscillations. The value function is monotone with respect to $x$, while it is non-monotone with respect to $i$. This is due to the non-parametric estimation result of the Markov chain, which is more flexible than the parametric models. This flexibility has both strong and weak points. The flexibility means the ability to fit to a variety of flow regimes, while it suggests importance of collecting data for the model estimation. In our case, we could estimate the Markov chain using the public data, while it would not be always possible to obtain such data for ungauged rivers. Coupling another methodology, such as another nonparametric hydrological modeling or a remote sensing technology, would be an option (Chouaib et al., 2019; Zhang et al., 2019).

The optimal harvesting amount $z^{*}$ equals $\bar{z}$ especially when the population is large and the discharge is small. The contrasting estimated Markov chains of 2018 and 2019 are reflected in the computational results of the triplet $\Phi, z^{*}, \lambda^{*}$: in the phase space $\Omega = J \times D$, the sub-domain where harvesting the population is optimal is wider in 2019 than in 2018, especially for the low flow regimes. The computed optimal observation frequencies are qualitatively the same between 2018 and 2019. This is in accordance with the estimated Markov chain in **Figures 3** and **4** where we can infer persistence of the low flow regimes in both years. More intensive observations would be optimal for the relatively low flow regimes. Doing nothing at observation times with less intensive observations is optimal for all the regimes if the population is small. It is also remarkable to see that harvesting the population is completely not optimal for the relatively large flow regimes, meaning that setting a sufficiently large $P$ can prevent the harvesting at such regimes, as also implied in the exact solution derived in **Section 3**.

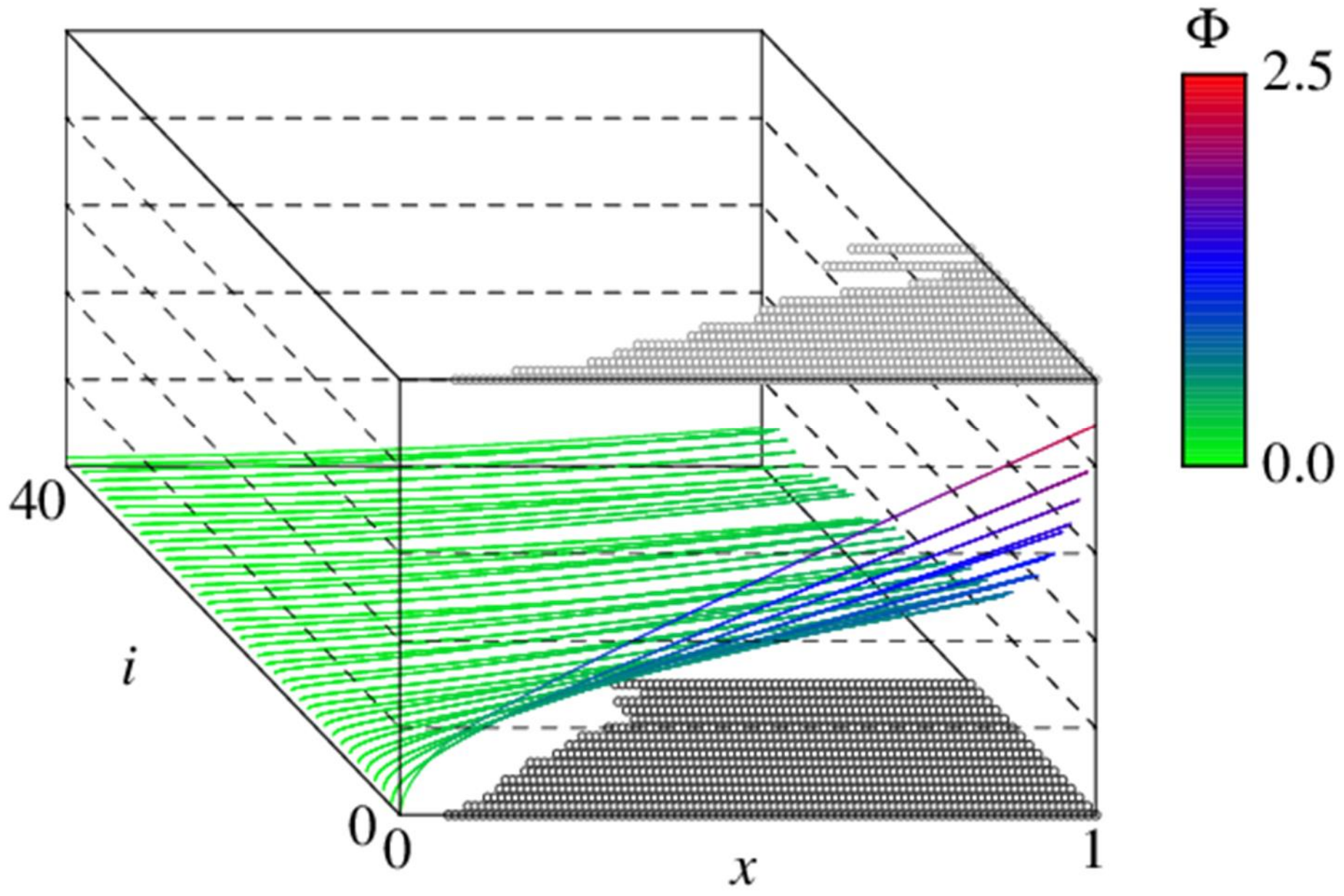

**Figure 5.** The value function $\Phi$ (solid curves), the optimal harvesting $z^*$ (black circles, lower surface), and the optimal observation intensity $\lambda^*$ (grey circles, upper surface) in $\Omega = J \times D$ using the data in 2018. Only the areas with $z^* = \bar{z}$ or $\lambda^* = \bar{\lambda}$ are plotted.

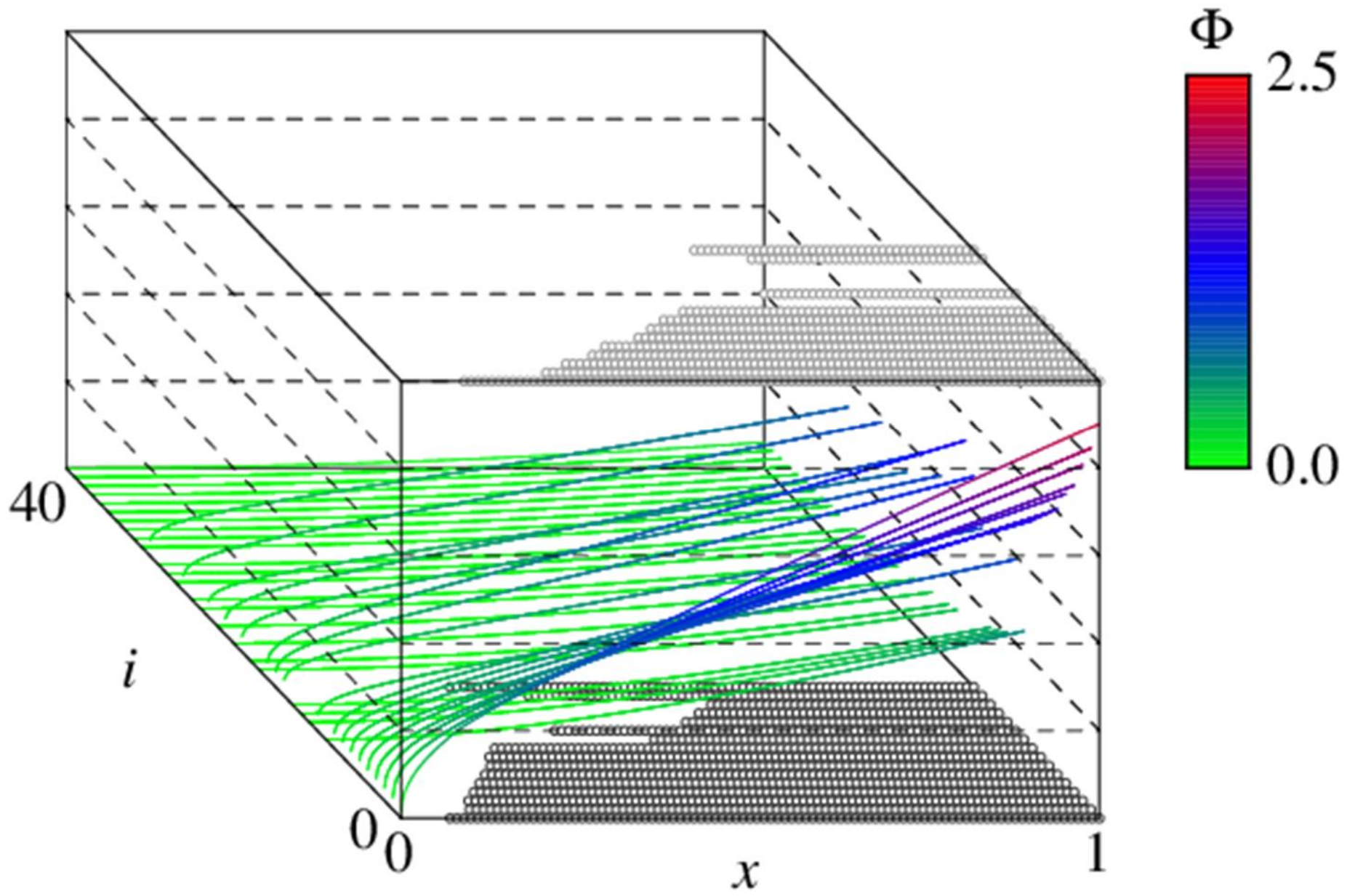

**Figure 6.** The value function $\Phi$, the optimal harvesting $z^*$, and the optimal observation intensity $\lambda^*$ in $\Omega = J \times D$ using the data in 2019. The same figure legends with **Figure 5**.

**Parameter sensitivity**

We present demonstrative computational examples on parameter dependence of the model. The parameter dependence is examined against the growth rate $\mu$ of the population, the discount rate $\delta$, and the power $m$. They are chosen because of their importance on the decision-making.

**Figures 7** through **9** show the computed optimal controls with different values of $\mu$, $\delta$, and $m$, respectively. The computational results suggest that the areas in the phase space $\Omega$ where more intensive observations and/or harvesting is optimal expand as the growth rate $\mu$ increases, the discount rate $\delta$ decreases, or $h$ becomes more concave (i.e., the power $m$ increases). Therefore, the results imply that longer-term river environmental management under a more rapid algae growth requires a more intensive monitoring and observations. Harvesting the population at each observation is optimal, except for a small portion of the phase space for the relatively low flow regimes with the linear and convex $h$. The concavity of $h$ means a sensitive reaction of the decision-maker on the small population; namely, the decision-maker with the concave utility should carefully observe and manage the algae population from the early growth stage of the population. This means that specifying the functional form of the disutility, which depends on the decision-makerȇs preference, is a key for designing the intervention policies. Such (dis-)utility functions can be specified through social science techniques (Pendleton and Mendelsohn, 2000; Swart and Zevenberg, 2018). At the current stage, it is important to construct a flexible model such that the optimal policy for a wide variety of the functional forms of the disutility can be computed. Our model can be such a candidate because of its potential ability to handle the convex, linear, and concave disutilities.

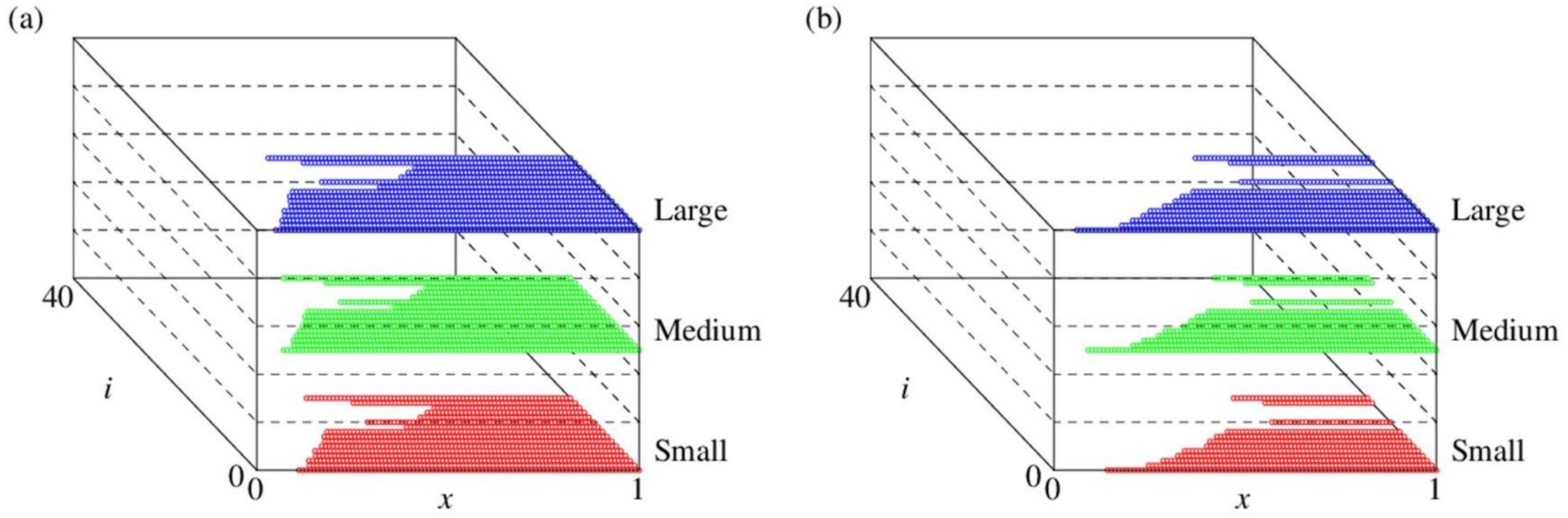

**Figure 7.** The optimal control variables (a) $z^*$ and (b) $\lambda^*$ for different values of $\mu$ (1/day): 0.35 (Small), 0.50 (Medium), and 0.65 (Large). Only the areas with $z^* = \overline{z}$ or $\lambda^* = \overline{\lambda}$ are plotted.

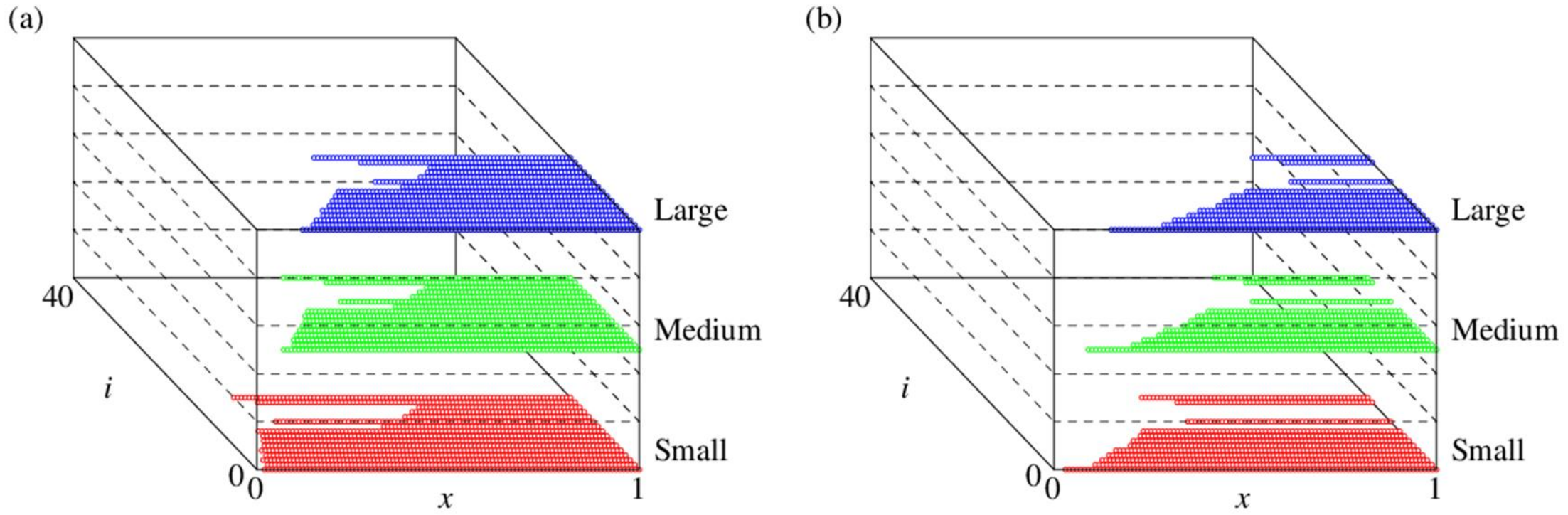

**Figure 8.** The optimal control variables (a) $z^*$ and (b) $\lambda^*$ for different values of $\delta$ (1/day): 0.10 (Small), 0.20 (Medium), and 0.30 (Large). Only the areas with $z^* = \overline{z}$ or $\lambda^* = \overline{\lambda}$ are plotted.

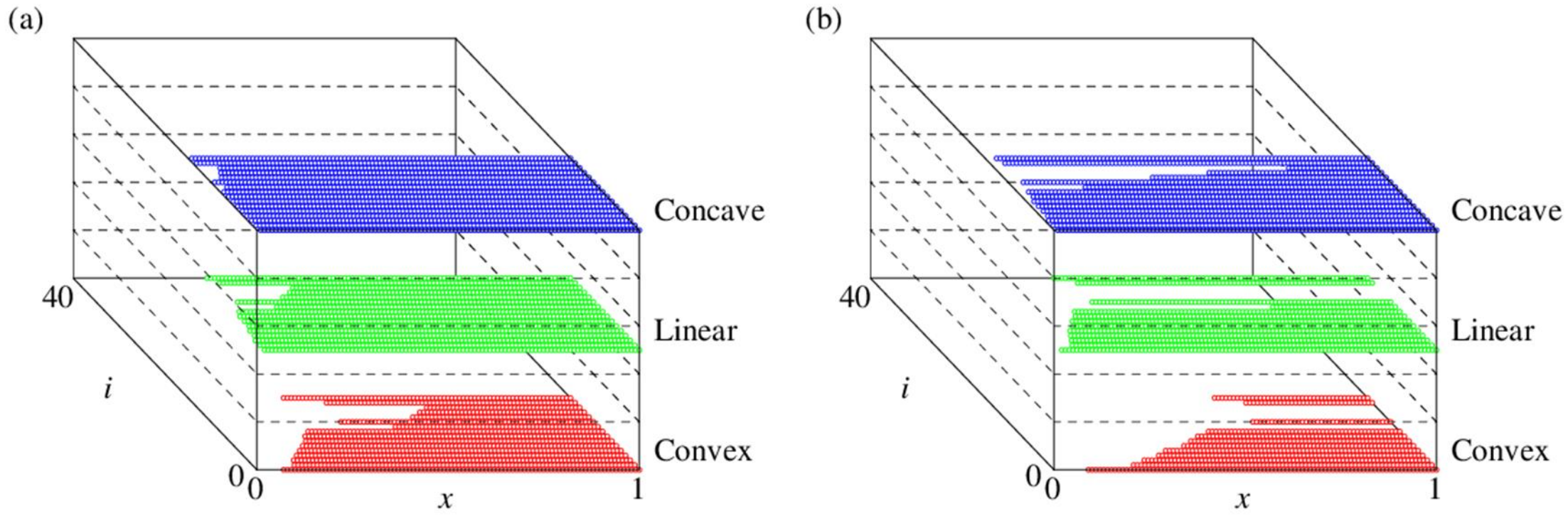

**Figure 9.** The optimal control variables (a) $z^*$ and (b) $\lambda^*$ for different values of $m$: 0.5 (Concave), 1.0 (Linear), and 2.0 (Convex). Only the areas with $z^* = \overline{z}$ or $\lambda^* = \overline{\lambda}$ are plotted.

**Comparison with the inflexible model.**

Finally, we consider a consequence of the different information structures coming from the flexibilities of the decision-maker. For this purpose, we compare the flexible and inflexible models; the latter is explained in **Appendix B**. Here, it is important to see that the optimality equation in the inflexible case has a similar mathematical form with that of the present flexible case (42) and (43). Therefore, the optimality equation in the inflexible case can also be handled with the finite difference scheme and an auxiliary time-advancing technique like (71) and (72). The two models contain the common coefficients and parameters. Hereafter, the super-scripts $(\mathrm{F})$ and $(\mathrm{I})$ represent the quantities for the flexible and inflexible cases, respectively. The value function in the inflexible case is computed using the same numerical scheme.

**Figure 10** plots the value function and the optimal controls for the inflexible case, showing that the harvesting and observation policies are less intensive than the flexible case. This would be due to the information structure of the inflexible case that the harvesting policy of the population at an observation time must follow the decision made at the previous observation time. Therefore, the decision-maker faces with a risk of river environmental management under relatively high flow regimes with $i > L$. This risk is mitigated in the flexible model because of making the decision based on the currently available information. Although not presented here, harvesting no population at all becomes optimal as the penalty parameter $P$ increases in the inflexible model, while it is not in the flexible model; in the latter case, the computed optimal policy under a large $P$ limit is almost the same with that in **Figure 6**.

By **Proposition C.1**, we have $\Phi^{(\mathrm{F})} \leq \Phi^{(\mathrm{I})}$ in $\Omega$. Set the Value of Information (VOI) as $V = \Phi^{(\mathrm{I})} - \Phi^{(\mathrm{F})} \geq 0$. The VOI defines the largest cost to be possibly paid to improve an inflexible observation scheme to a flexible observation scheme. This quantity can be readily computed from $\Phi^{(\mathrm{I})}$ and $\Phi^{(\mathrm{F})}$ by the definition. **Figure 11** plots the VOI $V = \Phi^{(\mathrm{I})} - \Phi^{(\mathrm{F})}$ against different values of the penalty parameter $P$ appearing in the harvesting cost $K$. The computational results show that the VOI is increasing with respect to the population $x$, suggesting the importance of collecting information even at a cost if the algae bloom is expected to be severe. In addition, we see that $V$ is increasing with respect to $P$, meaning that the information availability is more crucial for the decision-maker having a larger $P$. The analysis from the standpoint of information collection thus emphasizes that risk management of the cleaning up activities in the high flow regimes needs to be carefully considered.

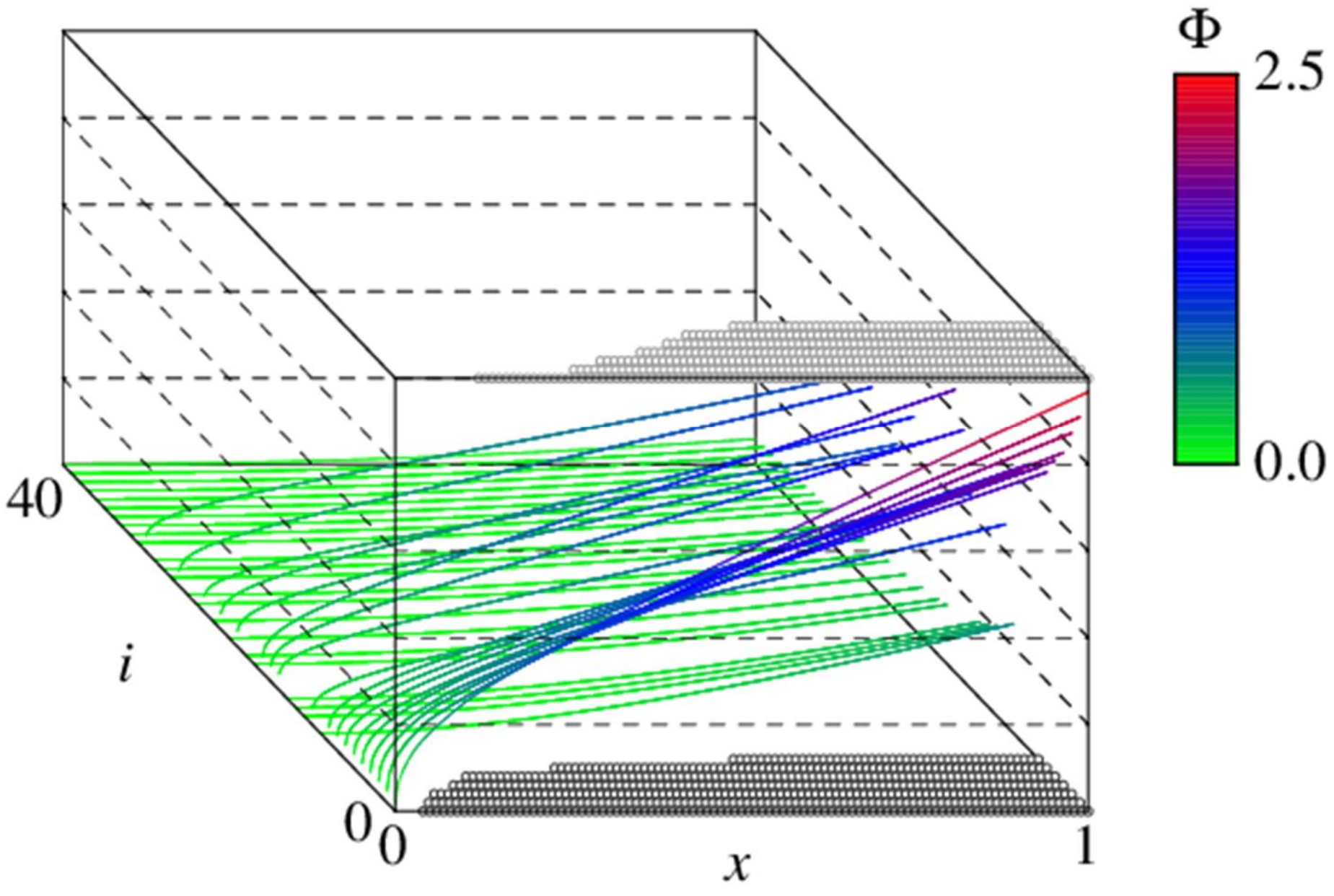

**Figure 10.** The value function $\Phi = \Phi^{(I)}$, the optimal harvesting $z^*$, and the optimal observation intensity $\lambda^*$ in $\Omega = J \times D$ using the data in 2019. The same figure legends with **Figure 5**.

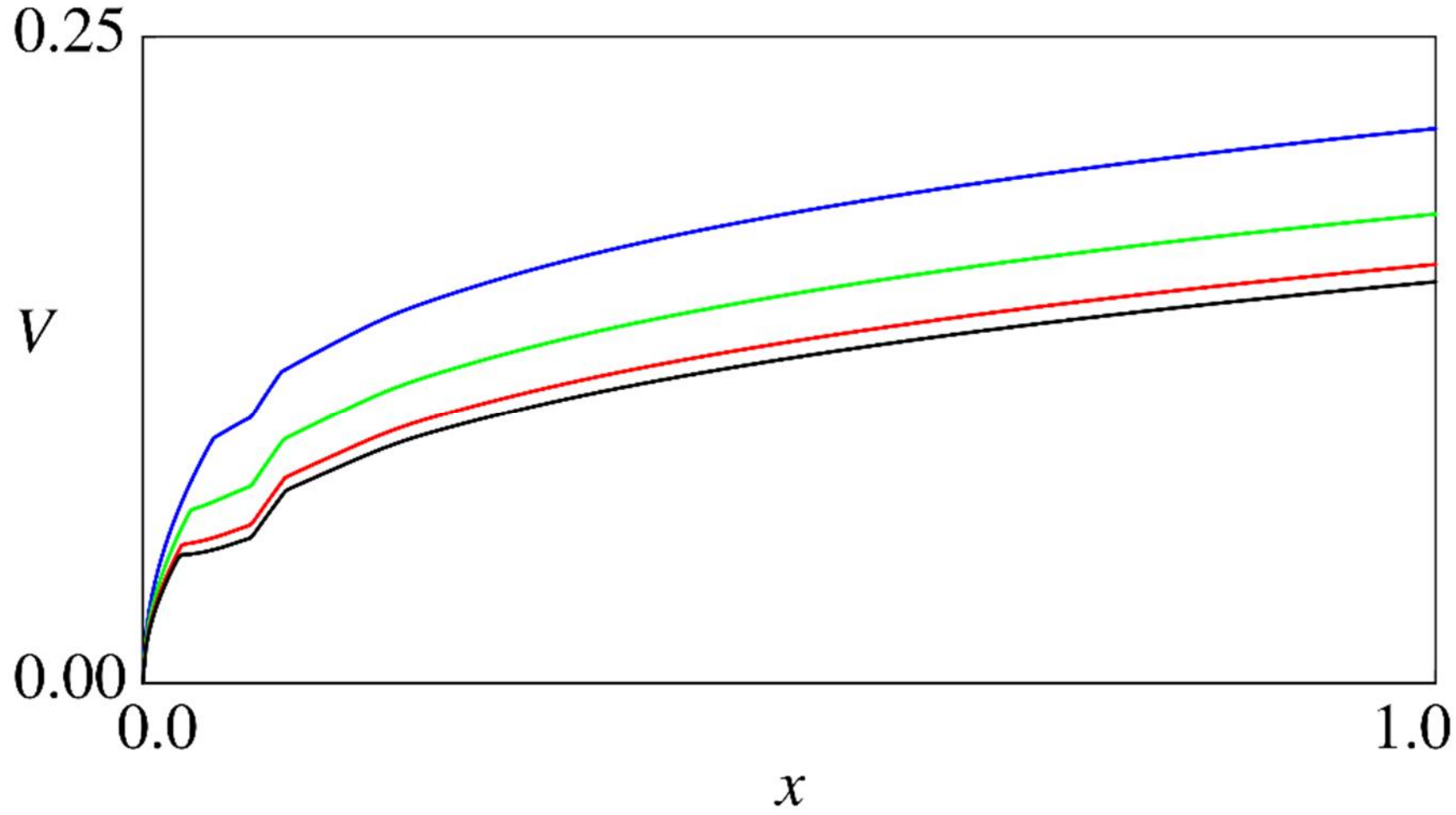

**Figure 11.** The VOI $V = \Phi^{(I)} - \Phi^{(F)}$ at the lowest flow regime $R_0$ for different values of $P$: $P = 5$ (Black), $P = 50$ (Red), $P = 200$ (Green), $P = 500$ (Blue).

## 5. Conclusions

We formulated a new stochastic control problem under discrete information, and mathematically and numerically analyzed its optimality equation and the optimal controls. The computational results based on the estimated parameter values suggested the optimal controls for the river environmental management depending on the preference of the decision-maker. The obtained results in this paper can be utilized for modeling and analysis of environmental management problems that many rivers around the world are facing with.

The presented mathematical framework is limited in the sense that only the single-species population dynamics are considered. Theoretically, it is possible to extend the presented dynamics to interacting multiple-species population dynamics representing aquatic ecosystems (Camara et al., 2019; Ghosh et al., 2019), but issues of computational costs may arise. Concerning the randomness of the observations, it may be possible to consider other distributions of the waiting time between each successive observation, like the gamma distribution. The optimal control problem with such a more flexible distribution may enable us to handle a wider range of situations, while at the same time we will lose the analytical tractability brought by the Poisson nature of the observation process. Considering more than one decision-makers, like fishery cooperatives and local governments, would lead to strategic game dynamics (Wang et al., 2016; Insley and Forsyth, 2020). In applications, there may be some implementation delay of the interventions (Zhang, 2018). Exploring environmental management policies under the above-mentioned situations can be interesting as well as important future topics.

**Acknowledgements**

JSPS research grant No. 18K01714, a grant for ecological survey of a life history of the landlocked ayu *Plecoglossus altivelis altivelis* from the Ministry of Land, Infrastructure, Transport and Tourism of Japan, and a research grant for young researchers in Shimane University support this research.

## Appendix A: Dynamic programming principle

In this appendix, we show the dynamic programming principle (39). Set the series of new admissible set of controls as $\mathcal{C}_n$ ( $n = 0,1,2,...$ ), such that $z_{n+j+1} = 0$ and $\lambda^{(n)} = \lambda^{(n+j)}$ ( $j = 0,1,2,...$ ). This is the set of controls adaptively choosing $\lambda^{(k)}$ only for $k \leq n$ and not harvesting for $\tau_n > t$. Clearly, $\mathcal{C}_n \subset \mathcal{C}$ is satisfied.

Set

$$\Phi_0(i,x) = \mathbb{E}^{i,x}\left[\int_0^\infty h(X_s)e^{-\delta s}\mathrm{d}s\right], \tag{74}$$

$$\Phi_n(i,x) = \inf_{u \in \mathcal{C}_n} \mathbb{E}^{i,x}\left[\int_0^\infty h(X_s)e^{-\delta s}\mathrm{d}s + \sum_{k=1}^{n} e^{-\delta \tau_k}\left(d + K\left(\alpha_{\tau_k}, X_{\tau_k}, z_k\right)\right)\right] \quad (n = 1,2,3,\ldots). \tag{75}$$

Notice that $\Phi_0$ is independent of the observation process.

We establish the dynamic programming principle in a step-by-step manner. At a first step, we show

***Lemma 1***

$$\lim_{n \to +\infty} \Phi_n = \Phi \quad \text{in} \quad J \times D. \tag{76}$$

**(Proof of Lemma 1)**

The inequality $\lim_{n \to +\infty} \Phi_n \geq \Phi$ is shown as follows. Choose a minimizer $u_n^* \in \mathcal{C}_n$ of $\Phi_n$ and its trivial extension in $\mathcal{C}$ as $\tilde{u}_n^* \in \mathcal{C}$ ($z_{n+j+1} = 0$ and $\lambda^{(n)} = \lambda^{(n+j)}$ ($j = 0,1,2,\ldots$)). We have

$$\begin{aligned}
\Phi(i,x) &\leq \phi\left(i,x;\tilde{u}_n^*\right) \\
&= \mathbb{E}^{i,x}\left[\int_0^\infty h(X_s)e^{-\delta s}\mathrm{d}s + \sum_{k=1}^{\infty} e^{-\delta \tau_k}\left(d + K\left(\alpha_{\tau_k}, X_{\tau_k}, z_k\right)\right)\right]_{u=\tilde{u}_n^*} \\
&= \mathbb{E}^{i,x}\left[\int_0^\infty h(X_s)e^{-\delta s}\mathrm{d}s + \sum_{k=1}^{n} e^{-\delta \tau_k}\left(d + K\left(\alpha_{\tau_k}, X_{\tau_k}, z_k\right)\right)\right]_{u=u_n^*} + d\mathbb{E}^{i,x}\left[\sum_{k=n+1}^{\infty} e^{-\delta \tau_k}\right]_{u=\tilde{u}_n^*} \\
&= \Phi_n(i,x) + d\mathbb{E}^{i,x}\left[\sum_{k=n+1}^{\infty} e^{-\delta \tau_k}\right]_{u=\tilde{u}_n^*}
\end{aligned} \tag{77}$$

The last term vanishes as $n \to \infty$, and we get the desired inequality by the Poisson observation.

We show $\lim_{n \to +\infty} \Phi_n \leq \Phi$. Set some $n \geq 1$ and a small $\varepsilon > 0$. Choose an $\varepsilon$-optimal control $u_\varepsilon \in \mathcal{C}$ and set its restriction in $\mathcal{C}_n$ as $u_{n,\varepsilon} \in \mathcal{C}_n$ (take first $n+1$ elements):

$$\begin{aligned}
\Phi(i,x) &\geq \phi\left(i,x;u_\varepsilon\right) - \varepsilon \\
&= \mathbb{E}^{i,x}\left[\int_0^\infty h(X_s)e^{-\delta s}\mathrm{d}s + \sum_{k=1}^{\infty} e^{-\delta \tau_k}\left(d + K\left(\alpha_{\tau_k}, X_{\tau_k}, z_k\right)\right)\right]_{u=u_\varepsilon} - \varepsilon \\
&\geq \mathbb{E}^{i,x}\left[\int_0^\infty h(X_s)e^{-\delta s}\mathrm{d}s + \sum_{k=1}^{n} e^{-\delta \tau_k}\left(d + K\left(\alpha_{\tau_k}, X_{\tau_k}, z_k\right)\right)\right]_{u=u_\varepsilon} - \varepsilon \\
&= \mathbb{E}^{i,x}\left[\int_0^\infty h(X_s)e^{-\delta s}\mathrm{d}s + \sum_{k=1}^{n} e^{-\delta \tau_k}\left(d + K\left(\alpha_{\tau_k}, X_{\tau_k}, z_k\right)\right)\right]_{u=u_{n,\varepsilon}} - \varepsilon
\end{aligned} \tag{78}$$

Then, we can just follow the argument of Lemma 10.1 of Øksendal and Sulem (2019). Notice that the stopping times $\{\tau_k\}_{k=1,2,3,\ldots}$ are not the control variable in our case, but are exogenously generated and satisfy a.s. $\tau_k = +\infty$ ($k \to +\infty$). In addition, the dynamics are confined in $J \times D$ without any constraints.

As the second step, we set another series of conditional expectations:

$$\varphi_0(i,x) = \mathbb{E}^{i,x}\left[\int_0^{\infty} h(X_s)e^{-\delta s}\mathrm{d}s\right], \tag{79}$$

$$\varphi_n(i,x) = \inf_{\lambda} \mathbb{E}^{i,x}\left[\int_0^{\tau_1} h(X_s)e^{-\delta s}\mathrm{d}s + e^{-\delta\tau_1}\mathcal{M}\varphi_{n-1}\left(\alpha_{\tau_1}, X_{\tau_1}\right)\right] \quad (n = 1,2,3,\ldots), \tag{80}$$

where

$$\mathcal{M}\varphi_{n-1}\left(\alpha_{\tau_1}, X_{\tau_1}\right) = \inf_{z_1}\left\{\varphi_{n-1}\left(\alpha_{\tau_1}, X_{\tau_1}\right) + d + K\left(\alpha_{\tau_1}, X_{\tau_1}, z_1\right)\right\}. \tag{81}$$

We show the following lemma.

***Lemma 2***

$$\Phi_n = \varphi_n \quad \text{in} \quad J \times D\,. \tag{82}$$

Combining **Lemmas 1** and **2** directly leads to

***Lemma 3***

$$\lim_{n\to+\infty}\varphi_n = \Phi \quad \text{in} \quad J \times D\,. \tag{83}$$

Then, **Lemma 3** leads to the dynamic programming principle.

**(Proof of Lemma 2)**

Firstly, we show $\varphi_n \le \Phi_n$. Set $n \ge 1$ and $u \in \mathcal{C}_n$. For $j \le n-1$, we have

$$\begin{aligned}\mathbb{E}^{i,x}\left[\varphi_{n-j}\left(\alpha_{\tau_j}, X_{\tau_j}\right)\right] &= \mathbb{E}^{i,x}\left[\int_{\tau_j}^{\tau_{j+1}} h(X_s)e^{-\delta(s-\tau_j)}\mathrm{d}s + e^{-\delta(\tau_{j+1}-\tau_j)}\mathcal{M}\varphi_{n-j-1}\left(\alpha_{\tau_{j+1}}, X_{\tau_{j+1}}\right)\right] \\ &\le \mathbb{E}^{i,x}\left[\begin{array}{l}\int_{\tau_j}^{\tau_{j+1}} h(X_s)e^{-\delta(s-\tau_j)}\mathrm{d}s \\ +e^{-\delta(\tau_{j+1}-\tau_j)}\left\{\varphi_{n-j-1}\left(\alpha_{\tau_{j+1}}, X_{\tau_{j+1}}\right) + d + K\left(\alpha_{\tau_{j+1}}, X_{\tau_{j+1}}, z_{j+1}\right)\right\}\end{array}\right].\end{aligned} \tag{84}$$

We can rewrite (84) as

$$\begin{aligned}&\mathbb{E}^{i,x}\left[e^{-\delta\tau_j}\varphi_{n-j}\left(\alpha_{\tau_j}, X_{\tau_j}\right) - e^{-\delta\tau_{j+1}}\varphi_{n-j-1}\left(\alpha_{\tau_{j+1}}, X_{\tau_{j+1}}\right)\right] \\ &\le \mathbb{E}^{i,x}\left[\int_{\tau_j}^{\tau_{j+1}} h(X_s)e^{-\delta s}\mathrm{d}s + e^{-\delta\tau_{j+1}}\left\{d + K\left(\alpha_{\tau_{j+1}}, X_{\tau_{j+1}}, z_{j+1}\right)\right\}\right].\end{aligned} \tag{85}$$

We then have the summation

$$
\begin{aligned}
&\sum_{j=0}^{n-1}\mathbb{E}^{i,x}\left[e^{-\delta\tau_j}\varphi_{n-j}\left(\alpha_{\tau_j},X_{\tau_j}\right)-e^{-\delta\tau_{j+1}}\varphi_{n-j-1}\left(\alpha_{\tau_{j+1}},X_{\tau_{j+1}}\right)\right]\\
&=\mathbb{E}^{i,x}\left[e^{-\delta\tau_0}\varphi_n\left(\alpha_{\tau_0},X_{\tau_0}\right)-e^{-\delta\tau_n}\varphi_0\left(\alpha_{\tau_n},X_{\tau_n}\right)\right]\\
&=\mathbb{E}^{i,x}\left[\varphi_n\left(i,x\right)-e^{-\delta\tau_n}\varphi_0\left(\alpha_{\tau_n},X_{\tau_n}\right)\right]\\
&=\varphi_n\left(i,x\right)-\mathbb{E}^{i,x}\left[e^{-\delta\tau_n}\varphi_0\left(\alpha_{\tau_n},X_{\tau_n}\right)\right]
\end{aligned}
\tag{86}
$$

and

$$
\begin{aligned}
&\sum_{j=0}^{n-1}\mathbb{E}^{i,x}\left[\int_{\tau_j}^{\tau_{j+1}}h\left(X_s\right)e^{-\delta s}\mathrm{d}s+e^{-\delta\tau_{j+1}}\left\{d+K\left(\alpha_{\tau_{j+1}},X_{\tau_{j+1}},z_{j+1}\right)\right\}\right]\\
&=\mathbb{E}^{i,x}\left[\int_0^{\tau_n}h\left(X_s\right)e^{-\delta s}\mathrm{d}s+\sum_{k=1}^{n}e^{-\delta\tau_k}\left(d+K\left(\alpha_{\tau_k},X_{\tau_k},z_k\right)\right)\right]
\end{aligned}
. \tag{87}
$$

In addition, by the strong Markov property, we get

$$
\begin{aligned}
\mathbb{E}^{i,x}\left[e^{-\delta\tau_n}\varphi_0\left(\alpha_{\tau_n},X_{\tau_n}\right)\right]&=\mathbb{E}^{i,x}\left[\mathbb{E}\left[e^{-\delta\tau_n}\varphi_0\left(\alpha_{\tau_n},X_{\tau_n}\right)\middle|\mathcal{F}_{\tau_n}\right]\right]\\
&=\mathbb{E}^{i,x}\left[\mathbb{E}^{\alpha_{\tau_n},X_{\tau_n}}\left[e^{-\delta\tau_n}\int_0^{\infty}h\left(X_s\right)e^{-\delta s}\mathrm{d}s\right]\right]\\
&=\mathbb{E}^{i,x}\left[\mathbb{E}^{\alpha_{\tau_n},X_{\tau_n}}\left[\int_0^{\infty}h\left(X_s\right)e^{-\delta\left(s+\tau_n\right)}\mathrm{d}s\right]\right]\\
&=\mathbb{E}^{i,x}\left[\mathbb{E}^{\alpha_{\tau_n},X_{\tau_n}}\left[\int_{\tau_n}^{\infty}h\left(X_{u-\tau_n}\right)e^{-\delta u}\mathrm{d}u\right]\right]\\
&=\mathbb{E}^{i,x}\left[\int_{\tau_n}^{\infty}h\left(X_s\right)e^{-\delta s}\mathrm{d}s\right]
\end{aligned}
\tag{88}
$$

Therefore, we arrive at

$$
\varphi_n\left(i,x\right)\le\mathbb{E}^{i,x}\left[\int_0^{\infty}h\left(X_s\right)e^{-\delta s}\mathrm{d}s+\sum_{k=1}^{n}e^{-\delta\tau_k}\left(d+K\left(\alpha_{\tau_k},X_{\tau_k},z_k\right)\right)\right]. \tag{89}
$$

Since $u\in\mathcal{C}_n$ is arbitrary, by taking, we obtain

$$
\varphi_n\left(i,x\right)\le\Phi_n\left(i,x\right). \tag{90}
$$

We show $\varphi_n\ge\Phi_n$. Set a small $\varepsilon>0$ and an $\varepsilon$-optimal control $u\in\mathcal{C}_n$:

$$
\mathbb{E}^{i,x}\left[\varphi_{n-j}\left(\alpha_{\tau_j},X_{\tau_j}\right)\right]\ge\mathbb{E}^{i,x}\left[\begin{aligned}&\int_{\tau_j}^{\tau_{j+1}}h\left(X_s\right)e^{-\delta\left(s-\tau_j\right)}\mathrm{d}s\\&+e^{-\delta\left(\tau_{j+1}-\tau_j\right)}\left\{\varphi_{n-j-1}\left(\alpha_{\tau_{j+1}},X_{\tau_{j+1}}\right)+d+K\left(\alpha_{\tau_{j+1}},X_{\tau_{j+1}},z_{j+1}\right)\right\}\end{aligned}\right]-\varepsilon. \tag{91}
$$

As in the discussion above, we get

$$
\varphi_n\left(i,x\right)\ge\mathbb{E}^{i,x}\left[\int_0^{\infty}h\left(X_s\right)e^{-\delta s}\mathrm{d}s+\sum_{k=1}^{n}e^{-\delta\tau_k}\left(d+K\left(\alpha_{\tau_k},X_{\tau_k},z_k\right)\right)\right]-n\varepsilon. \tag{92}
$$

Since $\varepsilon>0$ is arbitrary, we get the desired result by taking the limit $\varepsilon\to+0$ in (92).

Ĭ

## Appendix B: An inflexible model

We briefly analyze an inflexible model, which shares the same parameters and coefficient with the original flexible model, but with a different availability of information. A difference between the two models is the difference of the measurability of $z$: in the flexible model, we assume that $z_k$ is $\mathcal{F}_{\tau_k}$-measurable, while it is assumed to be only $\mathcal{F}_{\tau_{k-1}+}$-measurable in the inflexible model. Therefore, the decision-maker of the inflexible model decides the harvesting $z_k$ based on the information just after the previous observation time $t = \tau_{k-1}$. This means that the inflexible model assumes a less adaptive behavior of the decision-maker based on a partially feed-forward control. The optimality equation of the inflexible case is slightly different from that of the inflexible case.

The admissible set of control policies of the decision-maker in the inflexible model is denoted as $\mathcal{D}$, which is the set of the pair $u = (\lambda, z)$ of the process $\lambda = (\lambda_t)_{t \geq 0}$ and the sequence $z = \{z_k\}_{k=0,1,2,...}$: for each $k = 0,1,2,...$, both $\lambda^{(k)}$ and $z_k$ are $\mathcal{F}_{\tau_{k-1}+}$-measurable, $\lambda^{(k)} \in \Lambda = \{\bar{\lambda}, \underline{\lambda}\}$, and $z_k \in Z$. Clearly, $\mathcal{D}$ is non-empty and independent from $(\alpha_{0+}, X_{0+})$ as in the flexible model.

Following the partial observation framework (e.g., Pham and Tankov, 2009), the dynamic programming principle is heuristically given as

$$\Phi(i,x) = \inf_{\lambda_0, z} \mathbb{E}^{i,x}\left[\int_0^{\tau_1} h(\bar{X}_s) e^{-\delta s} \mathrm{d}s + \left\{\Phi\left(\alpha_{\tau_1}, (1-z)\bar{X}_{\tau_1}\right) + d + K\left(\alpha_{\tau_1}, \bar{X}_{\tau_1}, z\right)\right\} e^{-\delta \tau_1}\right], \quad (93)$$

where $\mathbb{E}^{i,x}$ represents the conditional expectation with respect to $\alpha$ and $N$, under the condition $(\alpha_{0+}, \bar{X}_{0+}) = (i,x) \in J \times D$. Here, the process $\bar{X} = (\bar{X}_t)_{t \geq 0}$ is identical to $X$ up to $t = \tau_1$, while no further intervention is carried out for $t > \tau_1$. The minimization is at most considered over the four cases: $(\lambda^{(0)}, z)$ is either $(\bar{\lambda}, 0)$, $(\bar{\lambda}, \bar{z})$, $(\underline{\lambda}, 0)$, or $(\underline{\lambda}, \bar{z})$. Therefore, the minimization in the right-hand side of (93) is carried out for the four discrete cases. In addition, because the observation and intervention are carried out only at the random and thus probabilistically distributed stopping time $\tau_1$, it is not possible to exchange the order between the minimization and the conditional expectation.

Introduce $\Omega = J \times D \times \Lambda \times Z$ and a function $\Psi : \Omega \to \mathbb{R}$, and set the notation

$$\Psi(i,x,r,y)=\mathbb{E}^{i,x}\left[\left.\begin{array}{l}\int_0^{\tau_1} h\left(\bar{X}_s\right)e^{-\delta s}\mathrm{d}s \\ +\left\{\Phi\left(\alpha_{\tau_1},(1-z)\bar{X}_{\tau_1}\right)+d+K\left(\alpha_{\tau_1},\bar{X}_{\tau_1},z\right)\right\}e^{-\delta\tau_1}\end{array}\right]\right|_{\lambda_0=r,z=y} \tag{94}$$
$$=\breve{\mathbb{E}}^{i,x}\left[\int_0^{\infty}\left\{h\left(\bar{X}_t\right)+r\Phi\left(\alpha_t,(1-y)\bar{X}_t\right)+d+K\left(\alpha_t,\bar{X}_t,y\right)\right\}e^{-(\delta+r)t}\mathrm{d}t\right]$$

and rewrite (93) as

$$\Phi(i,x)=\inf_{(r,y)\in\Lambda\times Z}\Psi(i,x,r,y). \tag{95}$$

By the Feynman-Kac formula (Zhu et al., 2015), we get the optimality equation:

$$\delta\Psi_i-f(i,x)\frac{\partial\Psi_i}{\partial x}+\sum_{\substack{j=0\\ j\neq i}}^{I}w_{i,j}\left(\Psi_i-\Psi_j\right)-h(x)+r\left[\Psi_i-\left\{\Phi_i\left((1-y)x\right)+d+K(i,x,y)\right\}\right]=0 \tag{96}$$

with (95), where $\Psi_i=\Psi_i(x,r,y)=\Psi(i,x,r,y)$. Notice the difference from the optimality equation of the flexible case (42)-(43). The last term in the inflexible case reflects the fact that the decision-making of the harvesting is not made at the latest observation.

We deduce the optimal control at each $\tau_k$ as

$$\left(\lambda^{(k),*},z_{k+1}^{*}\right)=\underset{(r,y)\in\Lambda\times Z}{\arg\min}\,\Psi\left(\alpha_{\tau_k},X_{\tau_k+},r,y\right). \tag{97}$$

The value function in the inflexible case shares the same continuity estimate with the flexible model. We define the viscosity solutions in the inflexible case, which is used in **Appendix C**.

***Definition B.1***

*A pair* $(\Psi,\Phi)$ *of continuous functions* $\Phi:J\times D\to\mathbb{R}$ *and* $\Psi:\Omega\to\mathbb{R}$ *is a viscosity solution to the optimality equation, the system (96) and (95), if :*

*(a) Sub-solution:* $\Phi(i,x)\leq\inf_{r,y}\Psi(i,x,r,y)$ *for all* $(i,x)\in J\times D$ *, and for each* $(i,r,y)\in J\times\Lambda\times Z$ *:*

$$\begin{aligned}&\delta\Psi(i,\bar{x},r,y)-f(i,\bar{x})\frac{\partial\varphi}{\partial x}(\bar{x})+\sum_{\substack{j=0\\ j\neq i}}^{I}w_{\bar{i},j}\left(\Psi(i,\bar{x},r,y)-\Psi(j,\bar{x},r,y)\right)-h(\bar{x})\\&+r\left[\Psi(i,\bar{x},r,y)-\left\{\Phi\left(i,(1-y)\bar{x}\right)+d+K(i,\bar{x},y)\right\}\right]\leq 0\end{aligned} \tag{98}$$

*for any continuously differentiable* $\varphi:D\to\mathbb{R}$ *and* $(\bar{i},\bar{x})\in J\times D$ *that is a local maximum of* $\Psi(\cdot,\cdot,r,y)-\varphi$.

*(b) Super-solution:* $\Phi(i,x)\geq\inf_{r,y}\Psi(i,x,r,y)$ *for all* $(i,x)\in J\times D$ *, and for each*

$(i,r,y) \in J \times \Lambda \times Z$ :

$$\begin{aligned} &\delta\Psi(i,\bar{x},r,y) - f(i,\bar{x})\frac{\partial \varphi}{\partial x}(\bar{x}) + \sum_{\substack{j=0 \\ j\neq i}}^{I} w_{\bar{i},j}\left(\Psi(i,\bar{x},r,y) - \Psi(j,\bar{x},r,y)\right) - h(\bar{x}) \\ &+r\left[\Psi(i,\bar{x},r,y) - \left\{\Phi(i,(1-y)\bar{x}) + d + K(i,\bar{x},y)\right\}\right] \geq 0 \end{aligned} \tag{99}$$

*for any continuously differentiable* $\varphi : D \to \mathbb{R}$ *and* $(\bar{i},\bar{x}) \in J \times D$ *that is a local minimum of* $\Psi(\cdot,\cdot,r,y) - \varphi$.

**Appendix C: Comparison of flexible and inflexible cases**

A byproduct of the comparison result of the value functions is the comparison between the value functions of the flexible and inflexible cases. We show that the value function of the flexible case is not larger than that of the inflexible case, demonstrating an advantage of doing the population management flexibly. The super-scripts $(\mathrm{F})$ and $(\mathrm{I})$ represent the quantities for the flexible and inflexible cases, respectively. We show the following proposition.

***Proposition C.1***

*We have* $\Phi^{(\mathrm{F})} \leq \Phi^{(\mathrm{I})}$ *in* $J \times D$.

**(Proof of Proposition C.1)**

It is sufficient to show that $\Phi^{(\mathrm{F})}$ with its associated auxiliary variable $\Psi^{(\mathrm{F})}$ is a viscosity sub-solution to the optimality equation in the inflexible case. Firstly, we have the sub-solution property of $\left(\Psi^{(\mathrm{F})},\Phi^{(\mathrm{F})}\right)$: $\Phi^{(\mathrm{F})}(i,x) \leq \inf_r \Psi^{(\mathrm{F})}(i,x,r)$ for all $(i,x) \in J \times D$, and for each $(i,r) \in J \times \Lambda$

$$\begin{aligned} &\delta\Psi^{(\mathrm{F})}(i,\bar{x},r) - f(i,\bar{x})\frac{\partial \varphi}{\partial x}(\bar{x}) + \sum_{\substack{j=0 \\ j\neq i}}^{I} w_{i,j}\left(\Psi^{(\mathrm{F})}(i,\bar{x},r) - \Psi^{(\mathrm{F})}(j,\bar{x},r)\right) - h(\bar{x}) \\ &+r\left[\Psi^{(\mathrm{F})}(i,\bar{x},r) - \mathcal{M}\Phi^{(\mathrm{F})}(i,\bar{x})\right] \leq 0 \end{aligned} \tag{100}$$

for any continuously differentiable $\varphi : D \to \mathbb{R}$ and $(\bar{i},\bar{x}) \in J \times D$ that is a local maximum of $\Psi^{(\mathrm{F})}(\cdot,\cdot,r) - \varphi$. For any $z \in Z$, we get

$$\mathcal{M}\Phi^{(\mathrm{F})}(i,\bar{x}) \leq \Phi^{(\mathrm{F})}(i,\bar{x}(1-z)) + K(i,\bar{x},z). \tag{101}$$

Combining (100) with (101) yields

$$\delta\Psi^{(\mathrm{F})}(i,\bar{x},r) - f(i,\bar{x})\frac{\partial\varphi}{\partial x}(\bar{x}) + \sum_{\substack{j=0\\ j\neq i}}^{I} w_{i,j}\left(\Psi^{(\mathrm{F})}(i,\bar{x},r) - \Psi^{(\mathrm{F})}(j,\bar{x},r)\right) - h(\bar{x})$$
$$+r\left[\Psi^{(\mathrm{F})}(i,\bar{x},r) - \left(\Phi^{(\mathrm{F})}(i,\bar{x}(1-z)) + K(i,\bar{x},z)\right)\right] \leq 0 \quad . \tag{102}$$

Secondly, we define the functions $\Upsilon^{(\mathrm{F})}(i,x,z) = \Phi^{(\mathrm{F})}(i,x)$ and $\Xi^{(\mathrm{F})}(i,x,r,z) = \Psi^{(\mathrm{F})}(i,x,r)$ for all $z \in Z$. The functions $\Upsilon^{(\mathrm{F})}$ and $\Xi^{(\mathrm{F})}$ are the value function and its auxiliary variable trivially augmented with the parameter $z \in Z$. The inequality (102) means that $\Xi^{(\mathrm{F})}$ satisfies the differential inequality required for viscosity sub-solutions in the inflexible case. Furthermore, by the viscosity property, we have $\Upsilon^{(\mathrm{F})}(i,x,z) \leq \inf_{r\in\Lambda} \Xi^{(\mathrm{F})}(i,x,r,z)$. Taking the infimum of the both sides of this inequality with respect to $z \in Z$ gives the inequality

$$\Phi^{(\mathrm{F})}(i,x) \leq \inf_{r\in\Lambda, z\in Z} \Xi^{(\mathrm{F})}(i,x,r,z). \tag{103}$$

Notice that we have $\Psi^{(\mathrm{F})}(i,x,r) = \inf_{z\in Z} \Xi^{(\mathrm{F})}(i,x,r,z)$ and $\Phi^{(\mathrm{F})}(i,x) = \inf_{z\in Z} \Upsilon^{(\mathrm{F})}(i,x,z)$ by their definitions. The proof is completed by (102) and (103).

Ĭ

**Declarations**

Funding

JSPS research grant No. 18K01714, a grant for ecological survey of a life history of the landlocked ayu Plecoglossus altivelis altivelis from the Ministry of Land, Infrastructure, Transport and Tourism of Japan, and a research grant for young researchers in Shimane University support this research.

Conflicts of interest/Competing interests (include appropriate disclosures)

Not applicable

Availability of data and material

Not applicable

Code availability

Not applicable

**References**

[1] Ang, T. K., Safuan, H. M., Sidhu, H. S., Jovanoski, Z., & Towers, I. N. (2019). Impact of Harvesting on a Bioeconomic Predatorï Prey Fishery Model Subject to Environmental Toxicant. Bulletin of Mathematical Biology, 1-20. https://doi.org/10.1007/s11538-019-00627-8

[2] Asatryan, V., Dallakyan, M., Vardanyan, T., Barseghyan, N., & Gabrielyan, B. (2019). An integrative evaluation of suitability of a river for natural reproduction of trout of Lake Sevan (Armenia). Environmental Science and Pollution Research, 1-10. https://doi.org/10.1007/s11356-019-07151-1

[3] Beville, S. T., Kerr, G. N., & Hughey, K. F. (2012). Valuing impacts of the invasive alga *Didymosphenia geminata* on recreational angling. Ecological Economics, 82, 1-10. https://doi.org/10.1016/j.ecolecon.2012.08.004

[4] Blubaugh, C. K., Hagler, J. R., Machtley, S. A., & Kaplan, I. (2016). Cover crops increase foraging activity of omnivorous predators in seed patches and facilitate weed biological control. Agriculture, Ecosystems & Environment, 231, 264-270. https://doi.org/10.1016/j.agee.2016.06.045

[5] Bolius, S., Wiedner, C., & Weithoff, G. (2019). Low invasion success of an invasive cyanobacterium in a chlorophyte dominated lake. Scientific Reports, 9(1), 8297. https://doi.org/10.1038/s41598-019-44737-8

[6] Bretschger, L., & Vinogradova, A. (2019). Best policy response to environmental shocks: Applying a stochastic framework. Journal of Environmental Economics and Management, 97, 23-41.

[7] Calder, J. (2018). Lecture notes on viscosity solutions. Available at http://www-users.math.umn.edu/~jwcalder/8590F18/viscosity_solutions.pdf.

[8] Camara, B. I., Yamapi, R., & Mokrani, H. (2019). Environmental stochastic effects on phytoplanktonï zooplankton dynamics. Nonlinear Dynamics, 96(3), 2013-2029. https://doi.org/10.1007/s11071-019-04902-0

[9] Castellano, R., & Cerqueti, R. (2012). Optimal consumption/investment problem with light stocks: A mixed continuous-discrete time approach. Applied Mathematics and Computation, 218(12), 6887-6898. https://doi.org/10.1016/j.amc.2011.12.065

[10] Chouaib, W., Alila, Y., & Caldwell, P. V. (2019). On the use of mean monthly runoff to predict the flowï duration curve in ungauged catchments. Hydrological Sciences Journal, 64(13), 1573-1587. https://doi.org/10.1080/02626667.2019.1657233

[11] Cormier, S. M., Zheng, L., Hayslip, G., & Flaherty, C. M. (2018). A field-based characterization of conductivity in areas of minimal alteration: A case example in the Cascades of northwestern United States. Science of The Total Environment, 633, 1657-1666. https://doi.org/10.1016/j.scitotenv.2018.02.018

[12] Crandall, M. G., Ishii, H., & Lions, P. L. (1992). Userôs guide to viscosity solutions of second order partial differential equations. Bulletin of the American Mathematical Society, 27(1), 1-67. ttps://doi.org/10.1090/S0273-0979-1992-00266-5

[13] Cullis, J. D., Gillis, C. A., Bothwell, M. L., Kilroy, C., Packman, A., & Hassan, M. (2012). A conceptual model for the blooming behavior and persistence of the benthic mat forming diatomDidymosphenia geminata in oligotrophic streams. Journal of Geophysical Research: Biogeosciences, 117(G2). 1-13. https://doi.org/10.1029/2011JG001891

[14] Di, M., Liu, X., Wang, W., & Wang, J. (2019). Manuscript prepared for submission to environmental toxicology and pharmacology pollution in drinking water source areas: Microplastics in the Danjiangkou Reservoir, China. Environmental Toxicology and Pharmacology, 65, 82-89. https://doi.org/10.1016/j.etap.2018.12.009

[15] Duwal, S., Winkelmann, S., Schütte, C., & von Kleist, M. (2015). Optimal treatment strategies in the context of ótreatment for preventionô against HIV-1 in resource-poor settings. PLoS Computational Biology, 11(4), e1004200. https://doi.org/10.1371/journal.pcbi.1004200

[16] Dyrssen, H., & Ekström, E. (2018). Sequential testing of a Wiener process with costly observations. Sequential Analysis, 37(1), 47-58. https://doi.org/10.1080/07474946.2018.1427973

[17] Federico, S., & Gassiat, P. (2014). Viscosity characterization of the value function of an investment-consumption problem in presence of an illiquid asset. Journal of Optimization Theory and Applications, 160(3), 966-991. https://doi.org/10.1007/s10957-013-0372-y

[18] Federico, S., Gassiat, P., & Gozzi, F. (2017). Impact of time illiquidity in a mixed market without full observation. Mathematical Finance, 27(2), 401-437. https://doi.org/10.1111/mafi.12101

[19] Gable, T. D., Windels, S. K., & Bruggink, J. G. (2017). The problems with pooling poop: confronting sampling method biases in wolf (*Canis lupus*) diet studies. Canadian Journal of Zoology, 95(11), 843-851. https://doi.org/10.1139/cjz-2016-0308

[20] Gassiat, P., Gozzi, F., & Pham, H. (2014). Investment/consumption problem in illiquid markets with regime-switching. SIAM Journal on Control and Optimization, 52(3), 1761-1786. https://doi.org/10.1137/120876976

[21] Ghosh, P., Das, P., & Mukherjee, D. (2019). Persistence and stability of a seasonally perturbed three species stochastic model of salmonoid aquaculture. Differential Equations and Dynamical Systems, 27(4), 449-465. https://doi.org/10.1007/s12591-016-0283-0

[22] Gillis, C. A., Dugdale, S. J., & Bergeron, N. E. (2018). Effect of discharge and habitat type on the occurrence and severity of *Didymosphenia geminata* mats in the Restigouche River, eastern Canada. Ecohydrology, 11(5), e1959. 1-7. https://doi.org/10.1002/eco.1959

[23] Gladyshev, M. I., & Gubelit, Y. I. (2019). Green Tides: New Consequences of the Eutrophication of Natural Waters (Invited Review). Contemporary Problems of Ecology, 12(2), 109-125. https://doi.org/10.1134/S1995425519020057

[24] Grandits, P., Kovacevic, R. M., & Veliov, V. M. (2019). Optimal control and the value of information for a stochastic epidemiological SIS-model. Journal of Mathematical Analysis and Applications, 476(2), 665-695. https://doi.org/10.1016/j.jmaa.2019.04.005

[25] Hopper, J. V., Pratt, P. D., McCue, K. F., Pitcairn, M. J., Moran, P. J., & Madsen, J. D. (2017). Spatial and temporal variation of biological control agents associated with *Eichhornia crassipes* in the Sacramento-San Joaquin River Delta, California. Biological Control, 111, 13-22. https://doi.org/10.1016/j.biocontrol.2017.05.005

[26] Hu, G., & Tian, K. (2019). On hybrid stochastic population models with impulsive perturbations. Journal of Biological Dynamics, 13(1), 385-406. https://doi.org/10.1080/17513758.2019.1609607

[27] Huang, Y., Li, X., Zhu, S., & Wu, Z. (2018). Identification of two hatch date-dependent cohorts and variation in early growth rate of mud carp (*Cirrhinus molitorella*) in the Pearl River Delta, China. Marine and Freshwater Research, 69(8), 1231-1237. https://doi.org/10.1071/MF17278

[28] Ide, J. I., Takeda, I., Somura, H., Mori, Y., Sakuno, Y., Yone, Y., & Takahashi, E. (2019). Impacts of hydrological changes on nutrient transport from diffuse sources in a rural river basin, western Japan. Journal of Geophysical Research: Biogeosciences, 124(8), 2565-2581. https://doi.org/10.1029/2018JG004513

[29] Insley, M., & Forsyth, P. A. (2019). Climate games. Whoôs on first? Whatôs on second?. LôActualite Economique 95. Available at https://cs.uwaterloo.ca/~paforsyt/Who_Feb_2019.pdf

[30] Ismail, Z., & Salim, K. (2013). determination of critical factors in implementing river clean up projects: A Malaysian case study. CLEANï Soil, Air, Water, 41(1), 16-23. https://doi.org/10.1002/clen.201200254

[31] Jiang, G. S., & Peng, D. (2000). Weighted ENO schemes for Hamilton--Jacobi equations. SIAM Journal on Scientific Computing, 21(6), 2126-2143. https://doi.org/10.1137/S106482759732455X

[32] Jiang, K., Merrill, R., You, D., Pan, P., & Li, Z. (2019). Optimal control for transboundary pollution under ecological compensation: a stochastic differential game approach. Journal of Cleaner Production, 241, 118391. https://doi.org/10.1016/j.jclepro.2019.118391

[33] Joshi, Y. K., Tortajada, C., & Biswas, A. K. (2012). Cleaning of the Singapore River and Kallang Basin in Singapore: economic, social, and environmental dimensions. International Journal of Water Resources Development, 28(4), 647-658. https://doi.org/10.1080/07900627.2012.669034

[34] Jost, J. (2005). Postmodern Analysis. Springer, Berlin, Heidelberg.

[35] Katsoulakis, M. A. (1994). Viscosity solutions of second order fully nonlinear elliptic equations with state

constraints. Indiana University Mathematics Journal, 43(2), 493-519. https://www.jstor.org/stable/24898087
[36] Katz, S. B., Segura, C., & Warren, D. R. (2018). The influence of channel bed disturbance on benthic Chlorophyll a: A high resolution perspective. Geomorphology, 305, 141-153. https://doi.org/10.1016/j.geomorph.2017.11.010
[37] Kazama, S., & Watanabe, K. (2018). Estimation of periphyton dynamics in a temperate catchment using a distributed nutrient-runoff model. Ecological Modelling, 367, 1-9. https://doi.org/10.1016/j.ecolmodel.2017.11.006
[38] King, S. A. (2014). Hydrodynamic control of filamentous macroalgae in a sub-tropical spring-fed river in Florida, USA. Hydrobiologia, 734(1), 27-37. https://doi.org/10.1007/s10750-014-1860-2
[39] Kvamsdal, S. F., Poudel, D., & Sandal, L. K. (2016). Harvesting in a fishery with stochastic growth and a mean-reverting price. Environmental and Resource Economics, 63(3), 643-663. https://doi.org/10.1007/s10640-014-9857-x
[40] Lempa, J. (2014). Bounded variation control of Itô diffusions with exogenously restricted intervention times. Advances in Applied Probability, 46(1), 102-120. https://doi.org/10.1239/aap/1396360105
[41] Lungu, E. M., & Øksendal, B. (1997). Optimal harvesting from a population in a stochastic crowded environment. Mathematical Biosciences, 145(1), 47-75. https://doi.org/10.1016/S0025-5564(97)00029-1
[42] Marten, A. L., & Moore, C. C. (2011). An options based bioeconomic model for biological and chemical control of invasive species. Ecological Economics, 70(11), 2050-2061. https://doi.org/10.1016/j.ecolecon.2011.05.022
[43] Oberman, A. M. (2006). Convergent difference schemes for degenerate elliptic and parabolic equations: Hamilton--Jacobi equations and free boundary problems. SIAM Journal on Numerical Analysis, 44(2), 879-895. https://doi.org/10.1137/S0036142903435235
[44] Øksendal, B., & Sulem, A. (2019). Applied Stochastic Control of Jump Diffusions. Springer, Cham.
[45] Pemy, M., Zhang, Q., & Yin, G. G. (2008). Liquidation of a large block of stock with regime switching. Mathematical Finance: An international Journal of Mathematics, Statistics and Financial Economics, 18(4), 629-648. https://doi.org/10.1111/j.1467-9965.2008.00351.x
[46] Pendleton, L., & Mendelsohn, R. (2000). Estimating recreation preferences using hedonic travel cost and random utility models. Environmental and Resource Economics, 17(1), 89-108. https://doi.org/10.1023/A:1008374423710
[47] Pham, H., & Tankov, P. (2009). A coupled system of integrodifferential equations arising in liquidity risk model. Applied Mathematics and Optimization, 59(2), 147-173. https://doi.org/10.1007/s00245-008-9046-9
[48] Pham, H., & Tankov, P. (2008). A model of optimal consumption under liquidity risk with random trading times. Mathematical Finance: An International Journal of Mathematics, Statistics and Financial Economics, 18(4), 613-627. https://doi.org/10.1111/j.1467-9965.2008.00350.x
[49] Picarelli, A., & Reisinger, C. (2020). Probabilistic error analysis for some approximation schemes to optimal control problems. Systems & Control Letters, 137, 104619. https://doi.org/10.1016/j.sysconle.2019.104619
[50] Pikosz, M., Messyasz, B., & GŃbka, M. (2017). Functional structure of algal mat (*Cladophora glomerata*) in a freshwater in western Poland. Ecological Indicators, 74, 1-9. https://doi.org/10.1016/j.ecolind.2016.09.041
[51] Pinheiro, S. (2018). Optimal harvesting for a logistic growth model with predation and a constant elasticity of variance. Annals of Operations Research, 260(1-2), 461-480. https://doi.org/10.1007/s10479-016-2242-0
[52] Reaver, N. G. F., Kaplan, D. A., Mattson, R. A., Carter, E., Sucsy, P. V., & Frazer, T. K. (2019). Hydrodynamic Controls on Primary Producer Communities in Spring Fed Rivers. Geophysical Research Letters, 46(9), 4715-4725. https://doi.org/10.1029/2019GL082571
[53] Ren, X., Yang, K., Che, Y., Wang, M., Zhou, L., & Chen, L. (2016). Spatial and temporal assessment of the initial

pattern of phytoplankton population in a newly built coastal reservoir. Frontiers of Earth Science, 10(3), 546-559. https://doi.org/10.1007/s11707-015-0543-2

[54] Rindorf, A., Dichmont, C. M., Thorson, J., Charles, A., Clausen, L. W., Degnbol, P., ... & Mace, P. (2017). Inclusion of ecological, economic, social, and institutional considerations when setting targets and limits for multispecies fisheries. ICES Journal of Marine Science, 74(2), 453-463. https://doi.org/10.1093/icesjms/fsw226

[55] Steffensen, K. D., & Mestl, G. E. (2016). Assessment of pallid sturgeon relative condition in the upper channelized Missouri River. Journal of Freshwater Ecology, 31(4), 583-595. https://doi.org/10.1080/02705060.2016.1196465

[56] Stratonovich, R. L. (2020). Theory of Information and its Value. Springer Nature, Switzerland AG.

[57] Strini, J. A., & Thonhauser, S. (2019). On a dividend problem with random funding. European Actuarial Journal, 9(2), 607-633. https://doi.org/10.1007/s13385-019-00208-y

[58] Swart, J. A., & Zevenberg, J. (2018). Utilitarian and nonutilitarian valuation of natural resources: a game-theoretical approach. Restoration Ecology, 26, S44-S53. https://doi.org/10.1111/rec.12504

[59] Tabayashi, Y., Miki, K., Godo, T., Yamamuro, M., & Kamiya, H. (2017). Multi-tracer identification of nutrient origin in the Hii River watershed, Japan. Landscape and Ecological Engineering, 13(1), 119-129. https://doi.org/10.1007/s11355-016-0307-5

[60] Tezuka, M., Ohgushi, K., Ogata, N., & Nagabayashi, H. (2014). Quantitative evaluation of slipperiness of attached algae for slip accident prevention in waterfront activities. Journal of Japan Society on Water Environment, 37(6), 239-250. https://doi.org/10.2965/jswe.37.239

[61] Tomlinson, L. M., Auer, M. T., Bootsma, H. A., & Owens, E. M. (2010). The Great Lakes *Cladophora* model: development, testing, and application to Lake Michigan. Journal of Great Lakes Research, 36(2), 287-297.https://doi.org/10.1016/j.jglr.2010.03.005

[62] Wang, H. (2001). Some control problems with random intervention times. Advances in Applied Probability, 33(2), 404-422. https://doi.org/10.1017/S0001867800010867

[63] Wang, S., Sun, P., & de Véricourt, F. (2016). Inducing environmental disclosures: A dynamic mechanism design approach. Operations Research, 64(2), 371-389. https://doi.org/10.1287/opre.2016.1476

[64] Wang, Y. (2018). Optimal Stopping with Discrete Costly Observations. U.U.D.M. Project Report 2018:33

[65] Winkelmann, S., Schütte, C., & Kleist, M. V. (2014). Markov control processes with rare state observation: theory and application to treatment scheduling in HIVï 1. Communications in Mathematical Sciences, 12(5), 859-877. https://dx.doi.org/10.4310/CMS.2014.v12.n5.a4

[66] Yaegashi, Y., Yoshioka, H., Unami, K., & Fujihara, M. (2018). A singular stochastic control model for sustainable population management of the fish-eating waterfowl *Phalacrocorax carbo*. Journal of Environmental Management, 219, 18-27. https://doi.org/10.1016/j.jenvman.2018.04.099

[67] Yang, B., Cai, Y., Wang, K., & Wang, W. (2019). Optimal harvesting policy of logistic population model in a randomly fluctuating environment. Physica A: Statistical Mechanics and its Applications, 526, 120817. https://doi.org/10.1016/j.physa.2019.04.053

[68] Yin, G. G., & Zhu, C. (2009). Hybrid Switching Diffusions: Properties and Applications. Springer, New York.

[69] Yoshioka, H. (2019). A simplified stochastic optimization model for logistic dynamics with control-dependent carrying capacity. Journal of Biological Dynamics, 13(1), 148-176. https://doi.org/10.1080/17513758.2019.1576927

[70] Yoshioka, H., & Tsujimura, M. (2020). Analysis and computation of an optimality equation arising in an impulse control problem with discrete and costly observations. Journal of Computational and Applied Mathematics, 366,

112399. https://doi.org/10.1016/j.cam.2019.112399

[71] Yoshioka, H., & Yaegashi, Y. (2018). Robust stochastic control modeling of dam discharge to suppress overgrowth of downstream harmful algae. Applied Stochastic Models in Business and Industry, 34(3), 338-354. https://doi.org/10.1002/asmb.2301

[72] Yoshioka H., Yaegashi Y., & Yoshioka Y. (2020). A discontinuous Hamiltonian approach for operating a dam-reservoir system in a river, ICIAE2020, March 26 to 30, 2020, Matsue, Kunibiki Messe. Proceedings, pp.226-233.

[73] Yoshioka H., & Yoshioka Y. (2020). Regime switching constrained viscosity solutions approach for controlling dam-reservoir systems. preprint. arXiv. http://arxiv.org/abs/2001.02070

[74] Yoshioka, H., Yoshioka, Y., Yaegashi, Y., Tanaka, T., Horinouchi, M., & Aranishi, F. (2019). Analysis and computation of a discrete costly observation model for growth estimation and management of biological resources. Computers & Mathematics with Applications. https://doi.org/10.1016/j.camwa.2019.08.017

[75] Zhang, Y., Chiew, F. H., Liu, C., Tang, Q., Xia, J., Tian, J., ... & Li, C. (2019). Can remotely sensed actual evapotranspiration facilitate hydrological prediction in ungauged regions without runoff calibration?. Water Resources Research, e24394. https://doi.org/10.1029/2019WR02623

[76] Zhang, Y. C. (2018). Entry-exit decisions with implementation delay under uncertainty. Applications of Mathematics, 63(4), 399-422. https://doi.org/10.21136/AM.2018.0205-17

[77] Zhu, C., & Yin, G. (2009). On hybrid competitive Lotkaï Volterra ecosystems. Nonlinear Analysis: Theory, Methods & Applications, 71(12), e1370-e1379. https://doi.org/10.1016/j.na.2009.01.166

[78] Zhu, C., Yin, G., & Baran, N. A. (2015). Feynmanï Kac formulas for regime-switching jump diffusions and their applications. Stochastics An International Journal of Probability and Stochastic Processes, 87(6), 1000-1032. https://doi.org/10.1080/17442508.2015.1019884